\documentclass[12pt, a4paper]{article}

\usepackage[latin1]{inputenc}
\usepackage{a4}

\usepackage{latexsym}
\usepackage{graphicx,color}
\usepackage{amsmath} 
\usepackage{amstext} 
\usepackage{amssymb}
\usepackage{mathrsfs} 
\usepackage{enumerate} 
\usepackage[numbers]{natbib} 

\newcommand{\grad}{\nabla}
\newcommand{\R}{\mathbb{R}} 
\newcommand{\N}{\mathbb{N}} 
\newcommand{\E}{\mathbb{E}} 
\newcommand{\F}{\mathcal{F}} 
\renewcommand{\P}{\mathbb{P}} 
\newcommand{\M}{\mathcal{M}} 
\renewcommand{\L}{\mathcal{L}}  
\newcommand{\spd}{\mathbb{S}}  
\newcommand{\I}{\mathcal{I}} 
\newcommand{\pb}{\mathcal{P}\mathcal{B}} 

\newtheorem{theorem}{Theorem}[section]

\newtheorem{corollary}[theorem]{Corollary}
\newtheorem{proposition}[theorem]{Proposition}
\newtheorem{lemma}[theorem]{Lemma}
\newtheorem{definition}[theorem]{Definition}
\newtheorem{example}{Example}[section]
\newtheorem{ass}{Assumption}[section]

\newtheorem{rem}{Remark}[section]

\begin{document}

\title{General existence and uniqueness of viscosity solutions for impulse control of jump-diffusions}
\author{Roland C. Seydel\footnote{Max Planck Institute for Mathematics in the Sciences, Inselstra\ss{}e 22, D-04103 Leipzig, Germany, \texttt{seydel@mis.mpg.de}. Support by the IMPRS ``Mathematics in the Sciences'' and the Klaus Tschira Foundation is gratefully acknowledged}}
\date{April 25, 2008\\This version: \today}
\maketitle

\begin{abstract}
General theorems for existence and uniqueness of viscosity solutions for Hamilton-Jacobi-Bellman quasi-variational inequalities (HJBQVI) with integral term are established. Such nonlinear partial integro-differential equations (PIDE) arise in the study of combined impulse and stochastic control for jump-diffusion processes. The HJBQVI consists of an HJB part (for stochastic control) combined with a nonlocal impulse intervention term.


Existence results are proved via stochastic means, whereas our uniqueness (comparison) results adapt techniques from viscosity solution theory. This paper is to our knowledge the first treating rigorously impulse control for jump-diffusion processes in a general viscosity solution framework; the jump part may have infinite activity. In the proofs, no prior continuity of the value function is assumed, quadratic costs are allowed, and elliptic and parabolic results are presented for solutions possibly unbounded at infinity.
\end{abstract}






\section{Introduction}

Impulse control is concerned with controlling optimally a stochastic process by giving it impulses at discrete intervention times. We consider a process $X$ evolving according to the controlled stochastic differential equation (SDE) 
\begin{equation}
\label{SDE}
dX_t
= \mu(t,X_{t-},\beta_{t-})\,dt + \sigma(t,X_{t-},\beta_{t-})\,dW_t + \int \ell(t,X_{t-},\beta_{t-},z)\, \overline{N}(dz,dt),
\end{equation}
for a standard Brownian motion $W$ and a compensated Poisson random measure $\overline{N}(dz,dt) = N(dz,dt) - 1_{|z| < 1} \nu(dz)dt)$ with (possibly unbounded) intensity measure $\nu$ (the jumps of a Lévy process), and the stochastic control process $\beta$ (with values in some compact set $B$). The impulses occur at stopping times $(\tau_i)_{i\ge 1}$, and have the effect 
\[
X_{\tau_{i}} = \Gamma(t,\check{X}_{\tau_{i}-}, \zeta_{i}),
\]
after which the process continues to evolve according to the controlled SDE until the next impulse. (Detailed notation and definitions are introduced in Section \ref{sec_model}.) We denote by $\gamma = (\tau_i,\zeta_i)_{i\ge 1}$ the impulse control strategy, and by $\alpha = (\beta,\gamma)$ a combined control consisting of a stochastic control $\beta$ and an impulse control $\gamma$. The aim is to maximise a certain functional, dependent on the impulse-controlled process $X^\alpha$ until the exit time $\tau$ (e.g., for a finite time horizon $T>0$, $\tau := \tau_S \wedge T$, where $\tau_S$ is the exit time of $X^\alpha$ from a possibly unbounded set $S$):
\begin{equation}
\tag{MAX} \label{MAX}
v(t,x) := \max_\alpha \E^{(t,x)} \left[ \int_t^\tau f(s,X^{\alpha}_s,\beta_s) ds + g(\tau, X^{\alpha}_{\tau})1_{\tau < \infty} + \sum_{\tau_j \le \tau} K(\tau_j,  \check{X}^{\alpha}_{\tau_j -}, \zeta_j ) \right]
\end{equation}
Here the negative function $K$ incorporates the impulse transaction costs, and the functions $f$ and $g$ are profit functions.\\

\textbf{Applications. } Typically, impulse control problems involve fixed (and variable) transaction costs, as opposed to singular control (only proportional transaction costs), and stochastic control (where there are no interventions). This is why an impulse control formulation is a natural framework for modelling incomplete markets in finance (the market frictions being the transaction costs). The standard reference for applications as well as for theory is certainly \citet{BeLi84}; as a more recent overview for jump-diffusions, \citet{OkSu05} can be helpful.

Although our results are quite general, we will mainly have financial applications in mind. Among those, we mention only the following (for some complements, see \cite{Ko99}):
\begin{itemize}
\item Option pricing with transaction costs \cite{Su01}, \cite{DaKo05}, \cite{Be06}
\item Tracking an optimal portfolio (\cite{PlSu04}, \cite{PaZa05}), or portfolio optimisation with transaction costs \cite{OkSu02}
\item Guaranteed Minimum Withdrawal Benefits (GMWB) in long-term insurance contracts \cite{ChFo07}
\item Control of an exchange rate by the Central Bank \cite{MuO98}, \cite{CaZa00}
\end{itemize}
The last application is a good example for combined control: there are two different means of intervention, namely interest rates (stochastic control) and FX market interventions. The stochastic control affects the process continuously (we neglect transaction costs here), and the impulses have fixed transaction costs, but have an immediate effect and thus can better react to jumps in the stochastic process.\\

\textbf{Quasi-variational inequality. } The standard approach for solution of (\ref{MAX}) is certainly to analyse the so-called Hamilton-Jacobi-Bellman quasi-variational inequality (HJBQVI) on the value function $v$, here in its parabolic form:
\begin{equation}
\tag{QVI} \label{QVI}
\min ( - \sup_{\beta \in B} \{ u_t + \L^\beta u + f^\beta \}, u - \M u ) = 0 \qquad \mbox{in } [0,T) \times S\\
\end{equation}
for $\L^\beta$ the infinitesimal generator of the SDE (\ref{SDE}) (where $y = (t,x)$), 
\begin{multline*}
\L^\beta u (y) = \frac{1}{2} tr \left(\sigma(y,\beta)\sigma^T(y,\beta)D_x^2 u(y)\right) + \langle \mu(y,\beta), \grad_x u(y) \rangle \\
+ \int u(t,x + \ell(y,\beta,z)) - u(y) - \langle \grad_x u(y), \ell(y,\beta,z) \rangle 1_{|z|< 1} \, \nu(dz),
\end{multline*}
and $\M$ the intervention operator selecting the momentarily best impulse,
\[
\M u(t,x) = \sup_\zeta \{ u(t,\Gamma(t,x,\zeta)) + K(t,x,\zeta) \}.
\]
In this paper, we will prove existence and uniqueness of viscosity solutions of (\ref{QVI}) (for a general introduction to viscosity solutions, we recommend the ``User's Guide'' by \citet{CIL}). Our goal is to establish a framework that can be readily used (and extended) in applications, without too many technical conditions. For other solution approaches to (\ref{QVI}) and impulse control problems, we refer to \citet{BeLi84}, \citet{Da93} and \citet{OkSu05}. 

(\ref{QVI}) is formally a nonlinear, nonlocal, possibly degenerate, second order parabolic partial integro-differential equation (PIDE). We point out that the investigated stochastic process is allowed to have jumps (jump-diffusion process), including so-called ``infinite-activity processes'' where the jump measure $\nu$ may be singular at the origin. (It can be argued that infinite-activity processes are a good model for stock prices, see, e.g., \citet{CT04}, \citet{EbKe95}.)\\

\textbf{Contribution. Literature. } Our main contribution is to rigorously treat general existence and uniqueness of viscosity solutions of impulse control QVI for the jump-diffusion process (\ref{SDE}). (An existence result can also be found in \citet{OkSu05}, but in the proof many technical details are omitted; the uniqueness part has not yet been covered to our best knowledge.) 
Apart from this main contribution, we investigate the QVI on the whole space $\R^d$ with appropriate boundary conditions (the boundary is in general not a null set of $\R^d$), and with only minimal assumptions on continuity of the value function (the continuity will then be a consequence of the viscosity solution uniqueness). Our solutions can be unbounded at infinity with arbitrary polynomial growth, provided appropriate conditions on the functions involved are satisfied, and superlinear transaction costs (e.g., quadratic) are allowed.


Let us now give a short overview on existing literature for QVI in the viscosity solution context. The book by \citet{OkSu05} is a recent overview on general stochastic control for jump processes, including QVI and viscosity solutions. The authors present rather general existence and uniqueness results for viscosity solutions of QVI. However, they assume that the value function is continuous (which is wrong in general), and they do not treat boundary values. Furthermore, their uniqueness proof (only in the diffusion case) assumes conditions that are not satisfied in a typical one-dimensional case. So the present paper can also serve as complement to their sketch of proof.
Other papers by the authors (mostly concerned with specific portfolio optimisation problems) are \cite{OkSu02}, \cite{AkMeSu96}, along with other references given in their book.

A complete and rigorous proof of existence and uniqueness in the viscosity solutions framework was given by K. Ishii \cite{Is93}, but only without jumps in the underlying stochastic process, and without stochastic control. Our approach to uniqueness is inspired by this paper. \citet{TaYo93} present an early proof for a finite time horizon (no exit time) including stochastic control and optimal switching, but again only for a continuous stochastic process and under rather restrictive assumptions. They use the continuity of the value function, the proof of which covers 11 pages. A more recent paper on impulse control (in a portfolio optimisation context) is \citet{LVMP07}. They too treat only a diffusion problem (without integral terms), but prove existence and uniqueness of discontinuous viscosity solutions. We extend their approach to existence to a general jump-diffusion context in the present paper.

As further references (without being exhaustive), we mention \citet{FlSo06} for viscosity solutions in general stochastic control, \citet{Ph98} for optimal stopping of controlled jump-diffusion processes, and \citet{Le89} for impulse control of piecewise deterministic processes. A very good account of discontinuous viscosity solutions and strong uniqueness for HJB equations is given in \citet{BaRo98}.

Because the Levy measure is allowed to have a singularity of second order at $0$, we cannot use the standard approach to uniqueness of viscosity solutions of PDE as used in \cite{Ph98} for optimal stopping, and in \cite{BeKaRe01} for singular control. For a more detailed discussion, we refer to \citet{JaKa06} (and the references therein), who were the first to propose a way to circumvent the problem for an HJB PIDE; see also the remark in the uniqueness section. 

For our proof of uniqueness, we will use and extend the framework as presented in the more recent paper \citet{BaIm07} (the formulation in \cite{JaKa06} does not permit an easy impulse control extension). The reader might also find helpful \citet{BaChIm08}.\\



\textbf{Contents. } The paper consists of 4 main sections. Section 2 presents the detailed problem formulation, the assumptions and a summary of the main result; some helpful properties of the investigated problem are derived in Section 3. The following Section 4 is concerned with existence of a QVI viscosity solution. After introducing the setting of our impulse control problem and several helpful results, we prove in Theorem \ref{theo_existence_qvi_parabolic} that the value function is a (discontinuous) QVI viscosity solution. The existence result for the elliptic QVI is deduced from the corresponding parabolic one. 
The last main section (Section 5) then starts with a reformulation of the QVI and several equivalent definitions for viscosity solutions. A maximum principle for impulse control is then derived, and used in a comparison result, which yields uniqueness and continuity of the QVI viscosity solution. The unboundedness of the domain is treated by a perturbation technique using strict viscosity sub- and supersolutions. The paper is complemented by a synthesis and summary at the end.

\textbf{Notation. } $\R^d$ for $d \ge 1$ is the Euclidean space equipped with the usual norm and the scalar product denoted by $\langle \cdot, \cdot \rangle$. For sets $A, B \subset \R^d$, the notation $A \subset \subset B$ (compactly embedded) means that $\overline{A} \subset B$, and $A^c = \R^d \setminus A$ is the complement of $A$. We denote the space of symmetric matrices $\subset \R^{d\times d}$ by $\spd^d$, $\ge$ is the usual ordering in $\R^{d \times d}$, i.e. $X \ge Y \Leftrightarrow X - Y$ positive semidefinite. $|\cdot|$ on $\spd^d$ is the usual eigenvalue norm. $C^2(\R^d)$ is the space of all functions two times continuously differentiable with values in $\R$, and as usual, $u_t$ denotes the time derivative of $u$. $L^2(\P;\R^d)$ is the Hilbert space of all $\P$-square-integrable measurable random variables with values in $\R^d$, the measure $\P^X = \P \circ X^{-1}$ is sometimes used to lighten notation.


\section{Model and main result}
\label{sec_model}

Let a filtered probability space $(\Omega,\F,(\F_t)_{t\ge 0},\P)$ satisfying the usual assumptions be given. Consider an adapted $m$-dimensional Brownian motion $W$, and an adapted independent $k$-dimensional pure-jump Lévy process represented by the compensated Poisson random measure $\overline{N}(dz,dt) = N(dz,dt) - 1_{|z| < 1} \nu(dz)dt)$, where as always $\int (|z|^2 \wedge 1)\, \nu(dz) < \infty$ for the Lévy measure $\nu$. 
We assume as usual that all processes are right-continuous. Assume the $d$-dimensional state process $X$ follows the stochastic differential equation with impulses
\begin{equation}
\label{viso_SDE}
\begin{split}
dX_t
&= \mu(t,X_{t-},\beta_{t-})\,dt + \sigma(t,X_{t-},\beta_{t-})\,dW_t + \int_{\R^k} \ell(t,X_{t-},\beta_{t-},z)\, \overline{N}(dz,dt),
\qquad \tau_i < t < \tau_{i+1}\\
X_{\tau_{i+1}}
&= \Gamma(t,\check{X}_{\tau_{i+1}-}, \zeta_{i+1}) \qquad \quad i \in \N_0
\end{split}
\end{equation}
for $\Gamma: \R^+_0 \times \R^{2d} \to \R^d$ measurable, and $\mu, \ell: \R^+_0 \times \R^d \times B \to \R^d$, $\sigma: \R^+_0 \times \R^d \times B \to \R^{d \times m}$ satisfying the necessary conditions such that existence and uniqueness of the SDE is guaranteed. $\beta$ is a càdlàg adapted stochastic control (where $\beta(t,\omega)\in B$, $B$ compact non-empty metric space), 
and $\gamma=(\tau_1,\tau_2, \ldots, \zeta_1, \zeta_2, \ldots)$ is the impulse control strategy, where $\tau_i$ are stopping times with $0 = \tau_0 \le \tau_1 \le \tau_2 \le \ldots$, and $\zeta_i$ are adapted impulses. The measurable transaction set $Z(t,x) \subset \R^d$ denotes the allowed impulses when at time $t$ in state $x$. We denote by $\alpha = (\beta,\gamma)$ the so-called combined stochastic control, where $\alpha \in A=A(t,x)$, the admissible region for the combined stochastic control. Admissible means here in particular 
that existence and uniqueness of the SDE be guaranteed, and that we only consider Markov controls (detailed in section \ref{sec_prob_props}). We further assume that all constant stochastic controls $\beta \in B$ are admissible.

The term $\check{X}^{\alpha}_{\tau_j -}$ denotes the value of the controlled $X^\alpha$ in $\tau_j$ with a possible jump of the stochastic process, but without the impulse, i.e., $\check{X}^{\alpha}_{\tau_j -} = X^{\alpha}_{\tau_j -} + \Delta X^{\alpha}_{\tau_j}$, where $\Delta$ denotes the jump of the stochastic process. So for the first impulse, this would be the process $\check{X}^{\alpha}_{\tau_1 -} = X^\beta_{\tau_1}$ only controlled by the continuous control. If two or more impulses happen to be at the same time (e.g., $\tau_{i+1} = \tau_i$), then (\ref{viso_SDE}) is to be understood as concatenation, e.g., $\Gamma(t,\Gamma(t,\check{X}_{\tau_{i}-}, \zeta_i),\zeta_{i+1})$. (The notation used here is borrowed from \citet{OkSu05}.)

The following conditions are sufficient for existence and uniqueness of the SDE for constant stochastic control $\beta$ (cf. \citet{GiSk72}, p. 273): There is an $x \in \R^d$ and $C>0$ with $\int | \ell(t,x,\beta,z) |^2 \nu(dz) \le C < \infty$  for all $t \in [0,T]$, $\beta \in B$. 
Furthermore, there exist $C>0$ and a positive function $b$ with $b(h)\downarrow 0$ for $h\to 0$ s.th.
\begin{equation} 
\label{eq_conds_ex_eind_SDE}
\begin{split}
|\mu(t,x,\beta) - \mu(&t,y,\beta)|^2 + |\sigma(t,x,\beta)) - \sigma(t,y,\beta)|^2 \\
& + \int | \ell(t,x,\beta,z) - \ell(t,y,\beta,z)|^2 \nu(dz)  \le  C |x-y|^2\\
|\mu(t+h,x,\beta) &- \mu(t,x,\beta)|^2 + |\sigma(t+h,x,\beta)) - \sigma(t,x,\beta)|^2 \\
 + \int | \ell(t&+h,x,\beta,z) - \ell(t,x,\beta,z)|^2 \nu(dz)  \le  C (1+|x|^2)b(h)
\end{split}
\end{equation}
These are essentially Lipschitz conditions, as one can easily check; see also \citet{Ph98}. A solution of (\ref{viso_SDE}) with non-random starting value will then have finite second moments, which is preserved after an impulse if for $X \in L^2(\P;\R^d)$, also $\Gamma(t,X,\zeta(t,X)) \in L^2(\P;\R^d)$. 
This is certainly the case if the impulses $\zeta_j$ are in a compact and $\Gamma$ continuous, which we will later assume.\\
We only assume that existence and uniqueness hold with constant stochastic control (by (\ref{eq_conds_ex_eind_SDE}) or weaker conditions as in \cite{GiSk72}), which guarantees that $A(t,x)$ is non-empty.

\begin{rem}
For existence and uniqueness of the SDE with arbitrary control process $\beta$, it is not sufficient to simply assume Lipschitz conditions on $x \mapsto \mu(t,x,\beta)$ etc. (even if the Lipschitz constants are independent of the control $\beta$), as done in \citet{Ph98}. If the control depends in a non-Lipschitzian or even discontinuous way on the current state, uniqueness or even existence of (\ref{viso_SDE}) might not hold. Compare also \citet{GiSk79}, p. 156.
\end{rem}

\begin{rem}
The version of the SDE in (\ref{viso_SDE}) is the most general form of all L\'evy SDE formulations currently used. If the support of $\nu$ is contained in the coordinate axes, then the one-dimensional components of the L\'evy process are independent, and the integral $\int_{\R^k} \ell(t,X_{t-},\beta_{t-},z)\, \overline{N}(dz,dt)$ can be splitted into several one-dimensional integrals (i.e. we obtain the form used in \citet{OkSu05}). If furthermore $\ell$ is linear in $z$, then we obtain the form as used in \citet{Pro04}, Th. V.32.
\end{rem}

The general (combined) impulse control problem is: find $\alpha =  (\beta, \gamma) \in A$ that maximizes the payoff starting in $t$ with $x$
\begin{equation}
J^{(\alpha)}(t,x) 
= \E^{(t,x)}  \left[ \int_t^{\tau} f(s,X^{\alpha}_s,\beta_s) ds + g(\tau, X^{\alpha}_{\tau})1_{\tau < \infty} + \sum_{\tau_j \le \tau} K(\tau_j, \check{X}^{\alpha}_{\tau_j -}, \zeta_j ) \right],
\end{equation}
where $f: \R^+_0 \times \R^d \times B \to \R$, $g: \R^+_0 \times \R^d \to \R$, $K: \R^+_0 \times \R^{2d} \to \R$ are measurable, and  $\tau = \tau_S = \inf\{ s \ge t: X^{\alpha}_s \not\in S \}$ is the exit time from some open set $S \subset \R^d$ (possibly infinite horizon), or $\tau = \tau_S \wedge T$ for some $T>0$ (finite horizon). Note that $X_s^\alpha$ is the value at $s$ after all impulses in $s$ have been applied; so ``intermediate values'' are not taken into account by this stopping time.\\
The value function $v$ is defined by 
\begin{equation}
v(t,x) = \sup_{\alpha \in A(t,x)} J^{(\alpha)}(t,x)
\end{equation}

We require the integrability condition on the negative parts of $f$, $g$, $K$
\begin{equation}
\label{eq_integrability}
\E^{(t,x)}  \left[ \int_t^{\tau} f^-(s,X^{\alpha}_s,\beta_s) ds + g^-(\tau, X^{\alpha}_{\tau})1_{\tau < \infty} + \sum_{\tau_j \le \tau} K^-(\tau_j, \check{X}^{\alpha}_{\tau_j -}, \zeta_j ) \right]
< \infty
\end{equation}
for all $\alpha \in A(t,x)$.\\

\textbf{Parabolic HJBQVI. } Let the impulse intervention operator $\M = \M^{(t,x)}$ be defined by
\begin{equation}
\M u(t,x) = \sup \{ u(t,\Gamma(t,x,\zeta)) + K(t,x,\zeta): \zeta \in Z(t,x) \}
\end{equation}
(define $\M u(t,x) = - \infty$ if $Z(t,x) = \emptyset$ -- we will exclude this case later on). 
Then the hope is to find the value function by investigating the following parabolic Hamilton-Jacobi-Bellman QVI ($T>0$ finite):
\begin{equation}
\label{qvi_parabolic}
\begin{split}
\min ( - \sup_{\beta \in B} \{ u_t + \L^\beta u + f^\beta \}, u - \M u ) &= 0 \qquad \mbox{in } [0,T) \times S\\
\min( u - g, u - \M u ) &= 0 \qquad \mbox{in } [0,T) \times (\R^d \setminus S)\\
u = \max(g,\M g, \M^2 g, &\ldots)  \qquad \mbox{on } \{T\} \times \R^d,
\end{split}
\end{equation}
for $\L^\beta$ the generator of $X$ in the SDE (\ref{viso_SDE}) for constant stochastic control $\beta$, and $f^\beta(\cdot) := f(\cdot,\beta)$. The generator $\L^\beta$ has the form ($y=(t,x)$):
\begin{multline}
\label{L_beta_parabolic}
\L^\beta u (y) = \frac{1}{2} tr \left(\sigma(y,\beta)\sigma^T(y,\beta)D_x^2 u(y)\right) + \langle \mu(y,\beta), \grad_x u(y) \rangle \\
+ \int u(t,x + \ell(y,\beta,z)) - u(y) - \langle \grad_x u(y), \ell(y,\beta,z) \rangle 1_{|z|< 1} \, \nu(dz),
\end{multline}

While the equation for $S$ in (\ref{qvi_parabolic}) can be motivated by Dynkin's formula and the fact that $v\ge \M v$ by the optimality of $v$, we have to argue why we consider the value function $v$ on $[0,T) \times \R^d$ instead of the interesting set $[0,T) \times S$: This is due to the jump term of the underlying stochastic process. While it is not possible to stay a positive time outside $S$ (we stop in $\tau_{S}$), it is well possible in our setting that the stochastic process jumps outside, but we return to $S$ by an impulse before the stopping time $\tau_S$ takes notice. Thus we must define $v$ outside $S$, to be able to decide whether a jump back to $S$ is worthwile. The boundary condition has its origin in the following necessary condition for the value function:
\begin{equation}
\label{boundary_inequality}
\min (v-g, v-\M v) = 0 \qquad\mbox{ in } [0,T) \times (\R^d \setminus S)
\end{equation}
This formalises that the controller can either do nothing (i.e., at the end of the day, the stopping time $\tau_{S}$ has passed, and the game is over), or can jump back into $S$, and the game continues. A similar condition holds at time $T<\infty$, with the difference that the controller is not allowed to jump back in time (as the device permitting this is not yet available to the public). So the necessary terminal condition can be put explicitly as
\begin{equation}
\label{terminal_condition}
v = \sup( g,\M g, \M^2 g , \ldots) \qquad\mbox{ on } \{T\} \times \R^d.\footnote{If $g$ is lower semicontinuous,  if $\M$ preserves this property and if the $\sup$ is finite, then it is well-known that $v(T,\cdot)$ is lower semicontinuous. Even if $g$ is continuous, $v(T,\cdot)$ need not be continuous}
\end{equation}

\begin{example} 
The impulse back from $\R^d \setminus S$ to $S$ could correspond to a capital injection into profitable business to avoid untimely default due to a sudden event.
\end{example}


\textbf{Well-definedness of QVI terms. } We need to establish conditions under which the terms $\L^\beta u$ and $\M u$ in the QVI (\ref{qvi_parabolic}) are well-defined. For $\M u$ compare the discussion of assumptions below.

The integral operator in (\ref{L_beta_parabolic}) is at the same time a differential operator of up to second order (if $\nu$ singular). This can be seen by Taylor expansion for $u \in C^{1,2}([0,T)\times \R^d)$ for some $0< \delta < 1$:
\begin{align}
&\int_{|z| < \delta} | u(t,x + \ell(x,\beta,z)) - u(t,x) - \langle \grad u(t,x), \ell(x,\beta,z) \rangle 1_{|z|< 1}| \, \nu(dz) \nonumber\\
\label{I1_Taylor_expansion} \le & \int_{|z| < \delta} |\ell(t,x,\beta,z)|^2 |D^2 u(t,\tilde{x})| \, \nu(dz)
\end{align}
for an $\tilde{x} \in B(x,\sup_{|z| < \delta} \ell(t,x,\beta,z))$. 

We take the Lévy measure $\nu$ as given. As for all Lévy measures, $\int (|z|^2 \wedge 1)\, \nu(dz) < \infty$. Let $p^* > 0$ be a number such that $\int_{|z|\ge 1} |z|^{p^*}\, \nu(dz) < \infty$, and let $q^* \ge 0$ be a number such that $\int_{|z| < 1} |z|^{q^*}\, \nu(dz) < \infty$ (think of $p^*$ as the \textit{largest} and $q^*$ the \textit{smallest} such number, even if it does not exist).
Then the expression $\sup_{\beta \in B} \{ u_t + \L^\beta u + f^\beta \}$ is well-defined, if, e.g., the following conditions are satisfied (depending on the singularity of $\nu$ in $0$ and its behaviour at infinity):
\begin{enumerate}
\item $u \in C^{1,2}([0,T)\times \R^d)$, $\sup_{\beta \in B} \sigma(t,x,\beta) < \infty$, $\sup_{\beta \in B} \mu(t,x,\beta) < \infty$, and $\sup_{\beta \in B} f(t,x,\beta) < \infty$ for all $(t,x) \in [0,T) \times S$
\item One of the following (all constants are independent of $\beta,z$, and the inequalities hold locally in $(t,x) \in [0,T) \times S$ if $C_{t,x}$ is used as constant\footnote{I.e., for all $(t,x) \in [0,T) \times S$, there exists a constant $C_{t,x}>0$ and a neighbourhood of $(t,x)$ in $[0,T) \times S$ such that $|\ell(t,x,\beta,z)| \le C_{t,x}$ for all $\beta,z$}):
  \begin{enumerate}[(a)]
  \item $|\ell(t,x,\beta,z)| \le C_{t,x}$ if $\nu(\R^d) < \infty$
  \item $|\ell(t,x,\beta,z)| \le C_{t,x}(|z|)$ and $|u(t,z)| \le C(1+|z|^{p^*})$ ($C$ independent of $t$)
  \item Or more generally, for $a,b>0$ such that $ab\le p^*$: $|\ell(t,x,\beta,z)| \le C_{t,x}(1+|z|^a)$, $|u(t,z)| \le C(1+|z|^b)$ on $|z|\ge 1$. Furthermore, $|\ell(t,x,\beta,z)| \le C_{t,x} |z|^{q^*/2}$ on $|z|< 1$
  \end{enumerate}
\end{enumerate}
In the following, we will assume one of the above conditions on $\ell$ (which is given), and define as in \citet{BaIm07} for a fixed polynomial function $R: \R^d \to \R$ satisfying condition 2.: 

\begin{definition}[Space of polynomially bounded functions]
$\pb = \pb([0,T)\times \R^d)$ is the space of all measurable functions $u: [0,T)\times \R^d \to \R$ such that
\[
|u(t,x)| \le C_u (1 + R(x))
\]
for a time-independent constant $C_u > 0$.
\end{definition}

As pointed out in \cite{BaIm07}, this function space $\pb$ is stable under lower and upper semicontinuous envelopes, and functions in $\pb$ are locally bounded. Furthermore, it is stable under (pointwise) limit operations, and the conditions for Lebesgue's Dominated Convergence Theorem are in general satisfied.\\

\textbf{Assumptions and main result. } Let us now formalise the conditions necessary for both the existence and the uniqueness proof in the following assumption (see also the discussion at the end of the section):
\begin{ass}
\label{ass_existence_uniqueness}
\begin{enumerate}[(V1)]
\item \label{Gamma_cont} $\Gamma$ and $K$ are continuous
\item \label{Transact_set_compact} The transaction set $Z(t,x)$ is non-empty and compact for each $(t,x) \in [0,T] \times \R^d$. For a converging sequence $(t_n,x_n)\to (t,x)$ in $[0,T]\times \R^d$ (with $Z(t_n,x_n)$ non-empty), $Z(t_n,x_n)$ converges to $Z(t,x)$ in the Hausdorff metric
\item \label{Ass_mu_etc_continuous} $\mu,\sigma,\ell,f$ are continuous in $(t,x,\beta)$ on $[0,T) \times \R^d \times B$
\item \label{Ass_ell_integrable} $\ell$ satisfies one of the conditions detailed above, and $\pb$ is fixed accordingly
\end{enumerate}
\end{ass}

Define $S_T := [0,T) \times S$ and its parabolic nonlocal ``boundary'' $\partial^+S_T := \left( [0,T) \times ( \R^d \setminus S ) \right) \cup \left( \{T\} \times \R^d \right)$. Further denote $\partial^*S_T := \left( [0,T) \times \partial S \right) \cup \left( \{T\} \times \overline{S} \right)$.

Apart from Assumption \ref{ass_existence_uniqueness}, we will need the following assumptions for the proof of existence:
\begin{ass}
\label{ass_theo_existence_qvi}
\begin{enumerate}[(E1)]
\item \label{ass_value_fct_in_C} The value function $v$ is in $\pb([0,T]\times\R^d)$ 
\item \label{ass_ex_g_cont} $g$ is continuous 
\item \label{ass_ex_v_cont_on_boundary} The value function $v$ satisfies for every $(t,x) \in \partial^* S_T$, and all sequences $(t_n,x_n)_n \subset [0,T) \times S$ converging to $(t,x)$:
\begin{gather*}
\liminf_{n\to\infty} v(t_n,x_n) \ge g(t,x)\\
\mbox{ If } v(t_n,x_n) > \M v(t_n,x_n) \;\forall n: \qquad \limsup_{n\to\infty} v(t_n,x_n) \le g(t,x)
\end{gather*}
\item \label{ass_ex_higher_variability} For all $\rho>0$, $(t,x) \in [0,T) \times S$, and sequences $[0,T) \times S \ni (t_n,x_n) \to (t,x)$, there is a constant $\hat{\beta}$ (not necessarily in $B$) and an $N \in \N$ such that for all $0<\varepsilon<1/N$ and $n \ge N$,
\[
\P ( \sup_{t_n \le s \le t_n + \varepsilon} | X^{\hat{\beta},t_n,x_n}_s - x_n | < \rho )
\le \P ( \sup_{t_n \le s \le t_n + \varepsilon} | X^{\beta^n,t_n,x_n}_s - x_n | < \rho ),
\]
where $X^{\beta^n,t_n,x_n}$ is the process according to SDE (\ref{viso_SDE}), started in $(t_n,x_n)$ and controlled by $\beta^n$.
\end{enumerate}
\end{ass}

The following assumptions are needed for the uniqueness proof ($\delta >0$):
\begin{ass}
\label{ass_uniqueness_qvi}
\begin{enumerate}[(U1)]
\item \label{U_ass_Udelta_ind_beta} If $\nu(\R^k)=\infty$: For all $(t,x) \in [0,T) \times S$, $U_\delta(t,x) := \{ \ell(x,\beta,z): |z| < \delta \}$ does not depend on $\beta$ 
\item \label{U_ass_Upositive} If $\nu(\R^k)=\infty$: For all $(t,x) \in [0,T) \times S$, $dist(x,\partial U_\delta(t,x))$ is strictly positive for all $\delta>0$ (or $\ell \equiv 0$)
\end{enumerate}
\end{ass}

\begin{ass}
\label{ass_comparison}
\begin{enumerate}[(B1)]
\item \label{ass_comparison_ints_bounded} $\int_{\R^d} |\ell(t,x,\beta,z) - \ell(t,y,\beta,z) |^2 \nu(dz) < C|x-y|^2$, $\int_{|z|\ge 1} |\ell(t,x,\beta,z) - \ell(t,y,\beta,z) | \nu(dz) < C|x-y|$, and all estimates hold locally in $t \in [0,T)$, $x,y$, uniformly in $\beta$
\item \label{ass_comparison_local_lipschitz} Let $\sigma(\cdot,\beta)$, $\mu(\cdot,\beta)$, $f(\cdot,\beta)$ be locally Lipschitz continuous, i.e. for each point $(t_0,x_0) \in [0,T) \times S$ there is a neigbourhood $U \ni (t_0,x_0)$ ($U$ open in $[0,T) \times \R^d$), and a constant $C$ (independent of $\beta$) such that $|\sigma(t,x,\beta) - \sigma(t,y,\beta)|\le C |x-y|$ $\forall (t,x),(t,y) \in U$, and likewise for $\mu$ and $f$
\end{enumerate}
\end{ass}

The space $\pb_p = \pb_p([0,T] \times \R^d)$ consists of all functions $u \in \pb$, for which there is a constant $C$ such that $|u(t,x)| \le C(1+|x|^p)$ for all $(t,x) \in [0,T] \times \R^d$.

Under the above assumptions, we can now formulate our main result in the parabolic case (the precise definition of viscosity solution is introduced in Section \ref{sec_existence}):

\begin{theorem}[Existence and uniqueness of a viscosity solution in the parabolic case]
\label{theo_existence_uniqueness_parabolic}
Let Assumptions \ref{ass_existence_uniqueness}, \ref{ass_theo_existence_qvi}, \ref{ass_uniqueness_qvi}, \ref{ass_comparison} be satisfied. Assume further that $v \in \pb_p([0,T] \times \R^d)$, and that there is a nonnegative $w \in \pb \cap C^2([0,T] \times \R^d)$ with $|w(t,x)|/|x|^p \to \infty$ for $|x| \to\infty$ (uniformly in $t$) and a constant $\kappa > 0$ such that
\begin{align*}
\min ( - \sup_{\beta \in B} \{ w_t + \L^\beta w + f^\beta \}, w - \mathcal{M}w ) &\ge \kappa  \qquad \mbox{in } S_T\\
\min( w - g, w - \M w ) &\ge \kappa \qquad \mbox{in } \partial^+ S_T.
\end{align*}
Then the value function $v$ is the unique viscosity solution of the parabolic QVI (\ref{qvi_parabolic}), and it is continuous on $[0,T] \times \R^d$.
\end{theorem}
For the proof of Theorem \ref{theo_existence_uniqueness_parabolic}, see Theorems \ref{theo_existence_qvi_parabolic} and \ref{viso_comparison_parabolic}.\\

\textbf{Elliptic HJBQVI. } For finite time horizon $T$, (\ref{qvi_parabolic}) is investigated on $[0,T] \times \R^d$ (parabolic problem). For infinite horizon, typically a discounting factor $e^{-\rho (t+s)}$ for $\rho > 0$ applied to all functions takes care of the well-definedness of the value function, e.g., $f(t,x,\beta) = e^{-\rho (t+s)} \tilde{f} (x,\beta)$. In this time-independent case, a transformation $u(t,x) = e^{-\rho (t+s)} w(x)$ gives us the elliptic HJBQVI to investigate
\begin{equation}
\label{qvi_elliptic}
\begin{split}
\min ( - \sup_{\beta \in B} \{ \L^\beta u + f^\beta \}, u - \M u ) &= 0 \qquad \mbox{in } S\\
\min( u - g, u - \M u ) &= 0 \qquad \mbox{in } \R^d \setminus S
\end{split}
\end{equation}
where the functions and variables have been appropriately renamed, and
\begin{align}
\L^\beta u (x) &= \frac{1}{2} tr \left(\sigma(x,\beta)\sigma^T(x,\beta)D^2 u(x)\right) + \langle \mu(x,\beta), \grad u(x) \rangle - \rho \,u(x) \nonumber \\
\label{L_beta_elliptic}
& \qquad + \int u(x + \ell(x,\beta,z)) - u(x) - \langle \grad u(x), \ell(x,\beta,z) \rangle 1_{|z|< 1} \, \nu(dz),\\
\M u(x) &= \sup \{ u(\Gamma(x,\zeta)) + K(x,\zeta): \zeta \in Z(x) \}.
\end{align}

Under the time-independent version of the assumptions above, an essentially identical existence and uniqueness result holds for the elliptic QVI (\ref{qvi_elliptic}). We refrain from repeating it, and instead refer to Sections \ref{subsec_existence_elliptic} and \ref{subsec_uniqueness_comparison} for a precise formulation.\\

\textbf{Discussion of the assumptions. } Of all assumptions, it is quite clear why we need the continuity assumptions, and they are easy to check.

By (V\ref{Ass_mu_etc_continuous}), (V\ref{Ass_ell_integrable}) and the compactness of the control set $B$, the Hamiltonian $\sup_{\beta \in B} u_t(t,x)  + \L^\beta u(t,x) + f^\beta(t,x)$ is well-defined and continuous in $(t,x) \in [0,T) \times S$ for $u \in \pb \cap C^{1,2}([0,T) \times \R^d)$. (This follows by sup manipulations, the (locally) uniform continuity and the DCT for the integral part.) 
Instead of (V\ref{Ass_mu_etc_continuous}), assuming the continuity of the Hamiltonian is sufficient for the existence proof. 
For the stochastic process $X_t$, condition (V\ref{Ass_ell_integrable}) essentially ensures the existence of moments.

By (V\ref{Gamma_cont}), (V\ref{Transact_set_compact}), we obtain that $\M u$ is locally bounded if $u$ locally bounded in $[0,T] \times \R^d$ (e.g., if $u \in \pb([0,T]\times \R^d)$). 
$\M u$ is even continuous if $u$ is continuous (so impulses preserve continuity properties), see Lemma \ref{lemma_Mu}.

In condition (V\ref{Transact_set_compact}), $Z(t,x) \not= \emptyset$ is necessary for the Hausdorff metric of sets to be well-defined, and to obtain general results on continuity of the value function (it is easy to construct examples of discontinuous value functions otherwise). The assumption is however no severe restriction, because we can set $Z(t,x)=\emptyset$ in the no-intervention region $\{ v > \M v \}$ without affecting the value function. The compactness of $Z(t,x)$ is not essential and can be relaxed in special cases --- this restriction is however of no practical importance.

Condition (E\ref{ass_ex_v_cont_on_boundary}) connects the combined control problem with the continuity of the stochastic control problem at the boundary. In this respect, Theorem \ref{theo_existence_uniqueness_parabolic}
roughly states that the value function is continuous except if there is a discontinuity on the boundary $\partial S$.
(E\ref{ass_ex_v_cont_on_boundary}) is typically satisfied if the stochastic process is regular at $\partial S$, as shown at the end of the section \ref{subsec_existence_parabolic}; see also \citet{FlSo06}, Theorem V.2.1, and the analytic approach in \citet{BaChIm08}. In particular, this condition excludes problems with true or \emph{de facto} state constraints, although the framework can be extended to cover state constraints.

(E\ref{ass_ex_higher_variability}) can be expected to hold because the control set $B$ is compact and the functions $\mu,\sigma,\ell$ are continuous in $(t,x,\beta)$ (V\ref{Ass_mu_etc_continuous}). The condition is very easy to check for a concrete problem --- it would be a lot more cumbersome to state a general result, especially for the jump part. 

\begin{example}
If $dX_t = \beta_t dt + dW_t$, and $\beta_t \in B = [-1,1]$, then $\hat{\beta} := 2$ is a possible choice for (E\ref{ass_ex_higher_variability}) to hold.
\end{example}

Assumption \ref{ass_uniqueness_qvi} collects some minor prerequisites that only need to be satisfied for small $\delta>0$ (see also the remark in the beginning of section \ref{subsec_uniqueness_prelims}), and the formulation can easily be adapted to a specific problem.

The local Lipschitz continuity in (B\ref{ass_comparison_ints_bounded}) and (B\ref{ass_comparison_local_lipschitz}) is a standard condition; (B\ref{ass_comparison_ints_bounded}) is satisfied if, e.g., the jump size of the stochastic process does not depend on $x$, or the conditions (\ref{eq_conds_ex_eind_SDE}) for existence and uniqueness of the SDE are satisfied for a constant $\beta$. Condition (B\ref{ass_comparison_ints_bounded}) can be relaxed if, e.g., $X$ has a state-dependent (finite) jump intensity -- the uniqueness proof adapts readily to this case. 

Certainly the most intriguing point is how to find a suitable function $w$ meeting all the requirements detailed in Theorem \ref{theo_existence_uniqueness_parabolic}. (This requirement essentially means that we have a strict supersolution.) We first consider the elliptic case of QVI (\ref{qvi_elliptic}). Here such a function $w$ for a $\kappa > 0$ can normally be constructed by $w(x) = w_1 | x |^{q} + w_2$ for suitable $w_i$ and $q>p$ (but still $w \in \pb$!). Main prerequisites are then 
\begin{enumerate}[(L1)]
\item \label{L_IR_positive} Positive interest rates: $\rho >\tilde{\kappa}$ for a suitably chosen constant $\tilde{\kappa}>0$ 
\item \label{L_FixedTAK} Fixed transaction costs: e.g., $K(x,\zeta)\le -k_0 < 0$
\end{enumerate}
If additionally we allow only impulses towards $0$, then $w - \M w > \kappa$ is easily achieved, as well as $w - g > \kappa$ (if we require that $g$ have a lower polynomial order than $w$). For a given bounded set, choosing $w_2$ large enough makes sure that $- \sup_{\beta \in B} \{ \L^\beta w + f^\beta \} > \kappa$ on this set (due to the continuity of the Hamiltonian and translation invariance in the integral). For $|x| \to \infty$, we need to impose conditions on $\tilde{\kappa}$ --- these depend heavily on the problem at hand, but can require the discounting factor to be rather large (e.g., for a geometric Brownian motion).

In the parabolic case, the same discussion holds accordingly, except that it is significantly easier to find a $w \in \pb([0,T] \times \R^d)$ satisfying assumption (L\ref{L_IR_positive}): By setting $w(t,x)=\exp(-\tilde{\kappa} t)( w_1 |x|^q + w_2 )$, we have $w_t = -\tilde{\kappa} w$ for arbitrarily large $\tilde{\kappa}$.


\section{Probabilistic properties}
\label{sec_prob_props}

In this section, we establish the Markov property of the impulse-controlled process, and derive a version of the dynamic programming principle (DPP) which we will need for the existence proof.

\textbf{Markov property. } For the Dynkin formula and several transformations, we need to establish the Markov property of the controlled process $X^\alpha$. To be more precise, we need to prove the strong Markov property of $Y^\alpha_t := (s + t, X^\alpha_{s+t})$ for some $s\ge 0$. First, by \citet{GiSk72} Theorem 1 in Part II.9, the Markov property for the uncontrolled process $X$ holds (for a similar result, see \citet{Pro04}, Th. V.32). If we consider the process $Y_t := (s+ t,X_{s+t})$, then the strong Markov property holds for $Y$.

We only consider Markov controls in the following, i.e. $\beta(t,\omega) = \beta(Y^\alpha(t-))$, the impulse times $\tau_i$ are exit times of $Y_t^\alpha$ (which makes sure they do not use past information), and $\zeta_i$ are $\sigma(\check{Y}^\alpha_{\tau_i-})$-measurable. (That Markov controls are sufficient is clear intuitively, because it is of no use to consider the past, if my objective function only depends on future actions and events, and my underlying process already has the Markov property.)

\begin{proposition}
Under the foregoing assumptions, the controlled process $Y^\alpha$ is a strong Markov process.
\end{proposition}

\textbf{Proof: } For a stopping time $T < \infty$ a.s., we have to show for all bounded measurable functions $h$ and for all $u\ge 0$ that $\E^y [h(Y^\alpha_{T+u}) | \F_T]$ is actually $Y^\alpha_T$-measurable. First note that without impulses, because $\beta(-t,\omega) = \beta(Y^\alpha(t-))$, the SDE solution $Y^\beta$ has the strong Markov property by the above cited results; we denote by $Y^\beta(y,t,t+u)$ this SDE solution started at $t\ge 0$ in $y$, evaluated at $t+u$.

Wlog, all $\tau_i \ge T$, e.g., $\tau_1$ first exit time after $T$. We split into the cases $A_0 := \{ \tau_1 > T + u \}$, $ A_1 := \{ \tau_1 \le T + u  < \tau_2 \}$, $A_2 := \{ \tau_2 \le T + u  < \tau_3 \}$, \ldots. The case $\tau_1 > T + u$ is clear by the above. For $A_1$, 
\begin{eqnarray*}
\E^y [1_{A_1} h(Y^\alpha_{T+u}) | \F_T] 
&=& \E^y [1_{A_1} \E^y [ h(Y^\alpha_{T+u})| \F_{\tau_1}  ] | \F_T] \\
&=& \E^y [1_{A_1} \E^y [ h(  Y^\beta(\Gamma(\check{Y}^\alpha_{\tau_1-},\zeta_1),\tau_1,T + u - \tau_1)) | \F_{\tau_1}  ] | \F_T].
\end{eqnarray*}
Because $\tau_1$ is first exit time after $T$ (and thus $T + u - \tau_1$ independent of $T$), and $Y^\alpha_{\tau_1-}$ includes the time information $\tau_1$ as first component, the SDE solution $Y^\beta(\Gamma(\check{Y}^\alpha_{\tau_1-},\zeta_1),\tau_1,T + u - \tau_1)$ depends only on $\check{Y}^\alpha_{\tau_1-}$. Thus we can conclude that there are measurable functions $g_1,\tilde{g_1}$ such that
\[
\E^y [1_{A_1} h(Y^\alpha_{T+u}) | \F_T] 
= \E^y [1_{A_1} g_1( \check{Y}^\alpha_{\tau_1-} ) | \F_T] 
= \E^y [1_{A_1} g_1( (\check{Y}^\alpha_{T+u-})^{\tau_1} ) | \F_T]
= \tilde{g_1} (Y^\alpha_T),
\]
where $(\check{Y}^\alpha_{T+u-})^{\tau_1}$ is the process stopped in $\tau_1$, which is a strong Markov process by \citet{Dy65}, Th. 10.2 and the fact that $\tau_1$ is the first exit time. 

For $A_2$, following exactly the same arguments, there are measurable functions $g_i,\tilde{g_i}$ such that
\begin{eqnarray*}
\E^y [1_{A_1} h(Y^\alpha_{T+u}) | \F_T] 
&=& \E^y [1_{A_1} g_2( \check{Y}^\alpha_{\tau_2-} ) | \F_T]\\
&=& \E^y [1_{A_1} \E^y [ g_2(  Y^\beta(\Gamma(\check{Y}^\alpha_{\tau_1-},\zeta_1),\tau_1,\tau_2 - \tau_1)) | \F_{\tau_1}  ] | \F_T]\\
&=& \E^y [1_{A_1} g_1( (\check{Y}^\alpha_{T+u-})^{\tau_1} ) | \F_T]
= \tilde{g_1} (Y^\alpha_T),
\end{eqnarray*}
where we have used that $\tau_2$ is the first exit time after $\tau_1$. The result follows by induction and the Dominated Convergence Theorem. \hfill$\Box$\\

\textbf{Dynamic Programming Principle. } An important insight into the structure of the problem is provided by Bellman's dynamic programming principle (DPP). Although the DPP is frequently used (see, e.g., \citet{OkSu05}, \citet{LVMP07}), we are only aware of the proof by \citet{TaYo93} in the impulse control case (for diffusions). We show here how the DPP can formally be derived from the Markov property.

By the strong Markov property of the controlled process, we have for a stopping time $\tilde{\tau}\le \tau$ ($\tau = \tau_S$ or $\tau = \tau_S \wedge T$):
\begin{eqnarray}
J^{(\alpha)}(t,x) 
&=& \E^{(t,x)}  \left\{ \int_t^{\tilde{\tau}} f(s,X^\alpha_s,\beta_s) ds 
+ \sum_{\tau_j < \tilde{\tau}} K(\tau_j,\check{X}^\alpha_{\tau_j -}, \zeta_j ) \right.\nonumber\\
&& \left. +\; \E^{(\tilde{\tau},\check{X}^\alpha_{\tilde{\tau}-})} \left[ 
\int_{\tilde{\tau}}^\tau f(s,X^\alpha_s,\beta_s) ds 
+ g(\tau, X^\alpha_{\tau}) 1_{\tau<\infty}
+ \sum_{\tilde{\tau} \le \tau_j \le \tau} K(\tau_j,\check{X}^\alpha_{\tau_j -}, \zeta_j ) \right]
\right\}\nonumber\\
& = & \label{dpp_ineq1} \E^{(t,x)}  \left\{ \int_t^{\tilde{\tau}} f(s,X^\alpha_s,\beta_s) ds 
+ \sum_{\tau_j < \tilde{\tau}} K(\tau_j,\check{X}^\alpha_{\tau_j -}, \zeta_j ) 
+ J^{(\alpha)}(\tilde{\tau},\check{X}^\alpha_{\tilde{\tau}-}) \right\}\\
& \le & \label{dpp_ineq2} \E^{(t,x)}  \left\{ \int_t^{\tilde{\tau}} f(s,X^\alpha_s,\beta_s) ds 
+ \sum_{\tau_j < \tilde{\tau}} K(\tau_j,\check{X}^\alpha_{\tau_j -}, \zeta_j ) 
+ v(\tilde{\tau},\check{X}^\alpha_{\tilde{\tau}-}) \right\}
\end{eqnarray}
Note especially that the second $J$ in (\ref{dpp_ineq1}) ``starts'' from $\check{X}_{\tilde{\tau}-}$, i.e. from $X$ before applying the possible impulses in $\tilde{\tau}$ -- this is to avoid counting a jump twice. $X_{\tilde{\tau}}$ instead of $\check{X}_{\tilde{\tau}-}$ in (\ref{dpp_ineq1}) would be incorrect (even if we replace the $=$ by a $\le$). However, $J^{(\alpha)}(\tilde{\tau},\check{X}_{\tilde{\tau}-}) \le v(\tilde{\tau},X_{\tilde{\tau}}) + K(\tilde{\tau},\check{X}_{\tilde{\tau}-},\zeta)1_{\{impulse\;in\;\tilde{\tau}\}}$ holds because a (possibly non-optimal) decision to give an impulse $\zeta$ in $\tilde{\tau}$ influences $J$ and $v$ in the same way. So we have the modified inequality
\begin{equation}
\label{dpp_inequality}
J^{(\alpha)}(t,x) 
\le \E^{(t,x)}  \left[ \int_t^{\tilde{\tau}} f(s,X^\alpha_s,\beta_s) ds 
+ \sum_{\tau_j \le \tilde{\tau}} K(\tau_j, \check{X}^\alpha_{\tau_j -}, \zeta_j ) 
+ v(\tilde{\tau},X^\alpha_{\tilde{\tau}}) \right].
\end{equation}
We will use both inequalities in the proof of Theorem \ref{theo_existence_qvi_parabolic} (viscosity existence).

The above considerations can be formalised in the well-known dynamic programming principle (DPP) (if the admissibility set $A(t,x)$ satisfies certain natural criteria): For all $\tilde{\tau}\le \tau$
\begin{equation}
\label{DPP}
v(t,x) 
= \sup_{\alpha \in A(t,x)} \E^{(t,x)}  \left[ \int_t^{\tilde{\tau}} f(s,X^\alpha_s,\beta_s) ds 
+ \sum_{\tau_j \le \tilde{\tau}} K(\tau_j, \check{X}^\alpha_{\tau_j -}, \zeta_j ) 
+ v(\tilde{\tau},X^\alpha_{\tilde{\tau}}) \right]
\end{equation}
The inequality $\le$ follows from (\ref{dpp_inequality}) and by approximation of the supremum. The inequality $\ge$ in (\ref{DPP}) is obvious from the definition of the value function, basically if two admissible strategies applied sequentially in time form a new admissible strategy (see also \citet{TaYo93}). (Of course a similar DPP can be derived from (\ref{dpp_ineq2}).) We do not use the DPP (\ref{DPP}) in our proofs.


\section{Existence}
\label{sec_existence}

In this section, we are going to prove the existence of a QVI viscosity solution in the elliptic and parabolic case. Because a typical impulse control formulation will include the time, we will first prove the existence for the parabolic form, then transforming the problem including time component into a time-independent elliptic one (the problem formulation permitting).


\subsection{Parabolic case}
\label{subsec_existence_parabolic}

Recall the definition of $S_T := [0,T) \times S$ and its parabolic ``boundary'' $\partial^+ S_T := \left( [0,T) \times ( \R^d \setminus S ) \right) \cup \left( \{T\} \times \R^d \right)$. We consider in this section the parabolic QVI in the form
\begin{equation}
\tag{\ref{qvi_parabolic}}
\begin{split}
\min ( - \sup_{\beta \in B} \{ u_t + \L^\beta u + f^\beta \}, u - \M u ) &= 0 \qquad \mbox{in } S_T\\
\min( u - g, u - \M u ) &= 0 \qquad \mbox{in } \partial^+ S_T 
\end{split}
\end{equation}
for $L^\beta$ from (\ref{L_beta_parabolic}) the integro-differential operator (or infinitesimal generator of the process $X$), and the intervention operator $\M$ selecting the optimal instantaneous impulse.

Let us now define what exactly we mean by a viscosity solution of (\ref{qvi_parabolic}). Let $LSC(\Omega)$ (resp., $USC(\Omega)$) denote the set of measurable functions on the set $\Omega$ that are lower semi-continuous (resp., upper semi-continuous). Let $T>0$, and let $u^*$ ($u_*$) define the upper (lower) semi-continuous envelope of a function $u$ on $[0,T] \times \R^d$, i.e. the limit superior (limit inferior) is taken only from within this set. Let us also recall the definition of $\pb$ encapsulating the growth condition from section \ref{sec_model}).

\begin{definition}[Viscosity solution]
A function $u \in \pb([0,T]\times\R^d)$ is a (viscosity) subsolution of (\ref{qvi_parabolic}) if for all $(t_0,x_0) \in [0,T]\times \R^d$ and $\varphi \in \pb \cap C^{1,2}([0,T)\times\R^d)$ with $\varphi(t_0,x_0) = u^*(t_0,x_0)$, $\varphi \ge u^*$ on $[0,T) \times \R^d$,
\begin{eqnarray*}
\min \left( - \sup_{\beta \in B} \left\{ \frac{\partial \varphi}{\partial t}  + \L^\beta \varphi + f^\beta \right\}, u^* - \mathcal{M}u^* \right) 
&\le & 0 \qquad \mbox{ in } (t_0,x_0) \in S_T\\
\min \left( u^* - g, u^* - \mathcal{M}u^* \right) 
&\le & 0 \qquad \mbox{ in } (t_0,x_0) \in \partial^+ S_T
\end{eqnarray*}
A function $u \in \pb([0,T]\times\R^d)$ is a (viscosity) supersolution of (\ref{qvi_parabolic}) if for all $(t_0,x_0) \in [0,T]\times \R^d$ and $\varphi \in \pb \cap C^{1,2}([0,T)\times\R^d)$ with $\varphi(t_0,x_0) = u_*(t_0,x_0)$, $\varphi \le u_*$ on $[0,T) \times \R^d$,
\begin{eqnarray*}
\min \left( - \sup_{\beta \in B} \left\{ \frac{\partial \varphi}{\partial t}  + \L^\beta \varphi + f^\beta \right\}, u_* - \mathcal{M}u_* \right) 
&\ge & 0 \qquad \mbox{ in } (t_0,x_0) \in S_T\\
\min \left( u_* - g, u_* - \mathcal{M}u_* \right) 
&\ge & 0 \qquad \mbox{ in } (t_0,x_0) \in \partial^+ S_T
\end{eqnarray*}
A function $u$ is a viscosity solution if it is sub- and supersolution.
\end{definition}

The conditions on the parabolic boundary are included inside the viscosity solution definition (sometimes called ``strong viscosity solution'', see, e.g., \citet{Is93}) because of the implicit form of this ``boundary condition''. In $T$, we chose the implicit form too, because otherwise the comparison result would not hold. The time derivative in $t=0$ is of course to be understood as a one-sided derivative.

Now we can state the main result of the section, the existence theorem:

\begin{theorem}[Viscosity solution: Existence]
\label{theo_existence_qvi_parabolic}
Let Assumptions \ref{ass_existence_uniqueness} and \ref{ass_theo_existence_qvi} be satisfied. Then the value function $v$ is a viscosity solution of (\ref{qvi_parabolic}) as defined above.
\end{theorem}

For the proof of Theorem \ref{theo_existence_qvi_parabolic}, we rely mainly on the proof given by \cite{LVMP07}, extending it to a general setting with jumps (compare also the sketch of proof in \citet{OkSu05}). We need a sequence of lemmas beforehand. The following lemma states first and foremost that the operator $\M$ preserves continuity. In a slightly different setting, the first two assertions can be found, e.g., in Lemma 5.5 of \cite{LVMP07}. 

\begin{lemma}
\label{lemma_Mu}
Let (V\ref{Gamma_cont}), (V\ref{Transact_set_compact}) be satisfied for all parts except (\ref{Mu_monotonous}). Let $u$ be a locally bounded function on $[0,T] \times \R^d$. Then
\begin{enumerate}[(i)]
\item \label{Mu_lsc} $\M u_* \in LSC([0,T] \times \R^d)$ and $\M u_* \le (\M u)_*$ 
\item \label{Mu_usc} $\M u^* \in USC([0,T]\times\R^d)$ and $(\M u)^* \le \M u^*$ 
\item \label{Mu_star_jump} If $u \le \M u$ on $[0,T] \times\R^d$, then $u^* \le \M u^*$ on $[0,T] \times\R^d$
\item \label{Mu_star_jump_m} For an approximating sequence $(t_n,x_n) \to (t,x)$, $(t_n,x_n) \subset [0,T]\times\R^d$ with $u(t_n,x_n) \to u^*(t,x)$: If $u^*(t,x) > \M u^*(t,x)$, then there exists $N \in \N$ such that $u(t_n,x_n) > \M u(t_n,x_n)$ $\forall n\ge N$
\item \label{Mu_monotonous} $\M$ is monotonous, i.e. for $u \ge w$, $\M u \ge \M w$. In particular, for the value function $v$, $\M^n v \le v$ for all $n\ge 1$
\end{enumerate}
\end{lemma}

\textbf{Proof: } (\ref{Mu_lsc}): Let $(t_n,x_n)_n$ be a sequence in $[0,T] \times \R^d$ converging to $(t,x)$. For an $\varepsilon>0$, select $\zeta^\varepsilon \in Z(t,x)$ with $u_*(t,\Gamma(t,x,\zeta^\varepsilon)) + K(t,x,\zeta^\varepsilon) + \varepsilon \ge \M u_*(t,x)$. Choose a sequence $\zeta_n \to \zeta^\varepsilon$ with $\zeta_n \in Z(t_n,x_n)$ for all $n$ (possible by Hausdorff convergence). Then
\begin{align*}
\liminf_{n\to\infty} \M u_*(t_n,x_n) 
&\ge \liminf_{n\to\infty} u_*(t,\Gamma(t_n,x_n,\zeta_n)) + K(t_n,x_n,\zeta_n)\\
&\ge u_*(t,\Gamma(t,x,\zeta^\varepsilon)) + K(t,x,\zeta^\varepsilon)
\ge \M u_*(t,x) - \varepsilon.
\end{align*}
The second assertion follows because $\M u \ge \M u_*$, and thus $(\M u)_* \ge (\M u_*)_* = \M u_*$.



(\ref{Mu_usc}): Fix some $(t,x) \in [0,T] \times \R^d$, and let $(t_n,x_n) \in [0,T] \times\R^d$ converge to $(t,x)$. Because of the upper semicontinuity of $u^*$ and continuity of $\Gamma$ and $K$, for each fixed $n$, the maximum in $\M u^*(t_n,x_n)$ is achieved, i.e., there is $\zeta_n \in Z(t_n,x_n)$ such that $\M u^*(t_n,x_n) = u^*(t_n, \Gamma(t_n,x_n,\zeta_n)) + K(t_n,x_n,\zeta_n)$. $(\zeta_n)_n$ is contained in a bounded set, thus has a convergent subsequence with limit $\hat{\zeta} \in Z(t,x)$ (assume that $\hat{\zeta} \not\in Z(t,x)$, then $dist(\hat{\zeta},Z(t,x))>0$, which contradicts the Hausdorff convergence).\\
We get
\begin{eqnarray*}
\M u^*(t,x) 
&\ge& u^*(t,\Gamma(t,x,\hat{\zeta})) + K(t,x,\hat{\zeta})
\ge \limsup_{n\to\infty} u^*(t_n,\Gamma(t_n,x_n,\zeta_n)) + K(t_n,x_n,\zeta_n)\\
&=& \limsup_{n\to\infty} \M u^*(t_n,x_n)
\end{eqnarray*}
The second assertion follows because $\M u \le \M u^*$, and thus $(\M u)^* \le (\M u^*)^* = \M u^*$.

(\ref{Mu_star_jump}): Follows immediately by (\ref{Mu_usc}): If $u\le \M u$, then $u^* \le (\M u)^* \le \M u^*$. 

(\ref{Mu_star_jump_m}): By contradiction: Assume $u(t_n,x_n) \le \M u(t_n,x_n)$ for infinitely many $n$. Then by convergence along a subsequence, \[ u^*(t,x) \le \limsup_{n\to\infty} \M u(t_n,x_n) \le (\M u)^*(t,x) \le \M u^*(t,x) .\]

(\ref{Mu_monotonous}): The monotonicity follows directly from the definition of $\M$. $\M v \le v$ is necessary for the value function $v$, because $v$ is already optimal, and $\M^n v \le v$ for all $n\ge 1$ then follows by induction.
\hfill$\Box$
\\

By (V\ref{Gamma_cont}), (V\ref{Transact_set_compact}) of section \ref{sec_model}, we obtain that $\M v(t,x) < \infty$ if $v$ locally bounded. This finiteness and the property that there is a convergent subsequence of $(\zeta_n)$ are sufficient for (\ref{Mu_usc}) (at least after reformulating the proof).

The existence proof frequently makes use of stopping times to ensure that a stochastic process $X$ (started at $x$) is contained in some (small) set. This works very well for continuous processes, because then for a stopping time $\tau = \inf \{ t\ge 0: |X_t - x| \ge \rho_1 \} \wedge \rho_2$, the process $|X_\tau - x|\le \rho_1$. For a process including (non-predictable) jumps however, $|X_\tau - x|$ may be $> \rho_1$. Luckily, L\'evy processes are stochastically continuous, which means that at least the \textit{probability} of $X_\tau$ being outside $\overline{B(x,\rho_1)}$ converges to $1$, if $\rho_2 \to 0$. Stochastic continuity as well holds for normal right-continuous Markov processes (see \citet{Dy65}, Lemma 3.2), and thus for our SDE solutions. 

The lemma (Lemma \ref{lemma_sup_cont} in the appendix) destined to overcome this problem essentially states the fact that stochastically continuous processes on a compact time interval are uniformly stochastically continuous. 
A further lemma (Lemma \ref{lemma_taumrho_not0} in the appendix) shows that for a continuously controlled process, stochastic continuity holds true uniformly in the control (this is of course a consequence of (E\ref{ass_ex_higher_variability})).

Now we are ready for the proof of the existence theorem. Recall that a necessary condition for the value function on the parabolic boundary $\partial^+ S_T$ is
\begin{equation}
\tag{\ref{boundary_inequality}}
\min (v-g, v-\M v) = 0 \qquad\mbox{ in } [0,T) \times (\R^d \setminus S).
\end{equation}


\textbf{Proof of Theorem \ref{theo_existence_qvi_parabolic}: } \textbf{$v$ is supersolution: } First, for any $(t_0,x_0) \in [0,T] \times \R^d$, the inequality $v(t_0,x_0) \ge \M v(t_0,x_0)$ holds, because otherwise an immediate jump would increase the value function. By Lemma \ref{lemma_Mu} (\ref{Mu_lsc}), $\M v_*(t_0,x_0) \le (\M v)_*(t_0,x_0) \le v_*(t_0,x_0)$. 

We then verify the condition on the parabolic boundary: Since we can decide to stop immediately, $v\ge g$ on $\partial^+ S_T$, so $v_* \ge g$ follows by the continuity of $g$ (outside $\overline{S}$) and requirement (E\ref{ass_ex_v_cont_on_boundary}) (if $x_0 \in \partial S$ or $t_0=T$).

So it remains to show the other part of the inequality
\begin{equation}
\label{supso_0}
- \sup_{\beta \in B} \left\{ \frac{\partial \varphi}{\partial t}  + \L^\beta \varphi + f^\beta \right\} \ge 0
\end{equation}
in a fixed $(t_0,x_0) \in [0,T)\times S$, for $\varphi \in \pb \cap C^{1,2}([0,T) \times \R^d)$, $\varphi(t_0,x_0) = v_*(t_0,x_0)$, $\varphi \le v_*$ on $[0,T) \times \R^d$.\\
From the definition of $v_*$, there exists a sequence $(t_n,x_n) \in [0,T) \times S$ such that $(t_n,x_n) \to (t_0,x_0)$, $v(t_n,x_n) \to v_*(t_0,x_0)$ for $n\to\infty$. By continuity of $\varphi$, $\delta_n := v(t_n,x_n) - \varphi(t_n,x_n)$ converges from above to $0$ as $n$ goes to infinity. Because $(t_0,x_0) \in [0,T) \times S$, there exists $\rho >0$ such that for $n$ large enough, $t_n < T$ and $B(x_n,\rho) \subset B(x_0, 2 \rho) = \{|y - x_0| < 2 \rho \} \subset S$.\\
Let us now consider the combined control with no impulses and a constant stochastic control $\beta \in B$, and the corresponding controlled stochastic process $X^{\beta, t_n,x_n}$ starting in $(t_n,x_n)$. Choose a strictly positive sequence $(h_n)$ such that $h_n \to 0$ and $\delta_n /h_n \to 0$ as $n\to\infty$. For 
\[ 
\bar{\tau}_n := \inf \{ s \ge t_n: | X^{\beta, t_n,x_n} - x_n |\ge \rho \} \wedge (t_n + h_n) \wedge T,
\]
we get by the strong Markov property and the Dynkin formula for $\rho$ sufficiently small ($\E^n = \E^{(t_n,x_n)}$ denotes the expectation when $X$ starts in $t_n$ with $x_n$):
\begin{eqnarray*}
v(t_n,x_n)
&\ge & \E^n \left[ \int_{t_n}^{\bar{\tau}_n} f(s,X^\beta_s, \beta) ds + v(\bar{\tau}_n, \check{X}^\beta_{\bar{\tau}_n-} )\right]\\
&\ge &  \E^n \left[ \int_{t_n}^{\bar{\tau}_n} f(s,X^\beta_s, \beta) ds + \varphi(\bar{\tau}_n, \check{X}^\beta_{\bar{\tau}_n-}) \right]\\
&=& \varphi(t_n,x_n) + \E^n \left[ \int_{t_n}^{\bar{\tau}_n} f(s,X^\beta_s, \beta) + \frac{\partial\varphi}{\partial t}(s,X^\beta_s) + \L^\beta \varphi(s,X^\beta_s)\; ds \right]
\end{eqnarray*}
Here, our assumptions on the SDE coefficients of $X$ were sufficient to apply Dynkin's formula because of the localizing stopping time $\bar{\tau}_n$. Using the definition of $\delta_n$, we obtain
\begin{equation}
\label{supso_1}
\frac{\delta_n}{h_n} 
\ge \E^n \left[ \frac{1}{h_n} \int_{t_n}^{\bar{\tau}_n} f(s,X^\beta_s, \beta) + \frac{\partial\varphi}{\partial t}(s,X^\beta_s) + \L^\beta \varphi(s,X^\beta_s)\; ds \right].
\end{equation}

Now, we want to let converge $n \to \infty$ in (\ref{supso_1}), but it is not possible to apply the mean value theorem because $s \mapsto f(s,X^\beta_s, \beta)$ (for fixed $\omega$) is in general not continuous. Select $\varepsilon \in (0,\rho)$. By Lemma \ref{lemma_sup_cont}, $\P(\sup_{t_n \le s \le r} | X^{\beta,t_n,x_n}_s - x_n | > \varepsilon) \to 0$ for $r \downarrow t_n$. Define now
\[
A_{n,\varepsilon} = \{ \omega: \sup_{t_n \le s \le t_n + h_n} | X_s^{\beta,t_n,x_n} - x_n | \le \varepsilon \}.
\]
and split the integral in (\ref{supso_1}) into two parts
\[
\left( \int_{t_n}^{t_n+h_n} \right) 1_{A_{n,\varepsilon}} + \left( \int_{t_n}^{\bar{\tau}_n} \right) 1_{A_{n,\varepsilon}^c}.
\]
On $A_{n,\varepsilon}$, for the integrand $G$ of the right hand side in (\ref{supso_1}),
\begin{eqnarray}
\left| G(t_0,x_0, \beta) - \frac{1}{h_n} \int_{t_n}^{t_n + h_n} G(s,X^\beta_s, \beta)\; ds \right| 
&\le& \frac{1}{h_n} \int_{t_n}^{t_n + h_n} | G(t_0,x_0, \beta) - G(s,X^\beta_s, \beta) |\; ds\nonumber\\
&\le& | G(t_0,x_0, \beta) - G(\hat{t}_{n,\varepsilon},\hat{x}_{n,\varepsilon}, \beta) |,
\end{eqnarray}
the latter because $G$ is continuous by (V\ref{Ass_mu_etc_continuous}) and assumption on $\varphi$, and the maximum distance of $| G(t_0,x_0, \beta) - G(\cdot,\cdot,\beta)|$ is assumed in a $(\hat{t}_{n,\varepsilon},\hat{x}_{n,\varepsilon}) \in [t_n,t_n+h_n] \times \overline{B(x_n,\varepsilon)}$.\\
On the complement of $A_{n,\varepsilon}$, 
\[
\frac{1}{h_n} \left( \int_{t_n}^{\bar{\tau}_n} \right) 1_{A_{n,\varepsilon}^c} 
\le \mbox{ess}\sup_{t_n \le s \le \bar{\tau}_n} \left| f(s,X^\beta_s, \beta) + \frac{\partial\varphi}{\partial t}(s,X^\beta_s) + \L^\beta \varphi(s,X^\beta_s) \right| 1_{A_{n,\varepsilon}^c},
\]
which is bounded by the same arguments as above and because a jump in $\bar{\tau}_n$ does not affect the essential supremum.

Because $h_n \to 0$ and $(t_n,x_n) \to (t_0,x_0)$ for $n\to\infty$ and by stochastic continuity, 
$\P(A_{n,\varepsilon})\to 1$, for all $\varepsilon>0$ or, equivalently, $1_{A_{n,\varepsilon}} \to 1$ almost surely. So by $n \to \infty$ and then $\varepsilon \to 0$, we can conclude by the dominated convergence theorem that $f(t_0,x_0, \beta) + \frac{\partial\varphi}{\partial t}(t_0,x_0) + \L^\beta \varphi(t_0,x_0) \le 0$ $\forall \beta \in B$, and thus (\ref{supso_0}) holds. \hfill$\Box$
\\


\textbf{$v$ is subsolution: } Let $(t_0,x_0) \in [0,T] \times \R^d$ and $\varphi \in \pb \cap C^{1,2}([0,T) \times \R^d)$ such that $v^*(t_0,x_0)=\varphi(t_0,x_0)$ and $\varphi \ge v^*$ on $[0,T) \times\R^d$. If $v^*(t_0,x_0)\le \M v^*(t_0,x_0)$, then the subsolution inequality holds trivially. So consider from now on the case $v^*(t_0,x_0) > \M v^*(t_0,x_0)$.

Consider $(t_0,x_0) \in \partial^+ S_T$. For an approximating sequence $(t_n,x_n)\to (t_0,x_0)$ in $[0,T] \times \R^d$ with $v(t_n,x_n)\to v^*(t_0,x_0)$, the relation $v(t_n,x_n)> \M v(t_n,x_n)$ holds by Lemma \ref{lemma_Mu} (\ref{Mu_star_jump_m}). Thus by the continuity of $g$ (outside $\overline{S}$) and requirement (E\ref{ass_ex_v_cont_on_boundary}) (if $x_0 \in \partial S$ or $t_0=T$), 
\[
g(t_0,x_0) = \lim_{n\to\infty} g(t_n,x_n) = \lim_{n\to\infty} v(t_n,x_n) = v^*(t_0,x_0)
\]

Now let us show that, for $v^*(t_0,x_0) > \M v^*(t_0,x_0)$,
\begin{equation}
\label{subso_0}
- \sup_{\beta \in B} \left\{ \frac{\partial \varphi}{\partial t}  + \L^\beta \varphi + f^\beta \right\} \le 0
\end{equation}
in $(t_0,x_0) \in [0,T) \times S$. We argue by contradiction and assume that there is an $\eta>0$ such that
\begin{equation}
\label{subso_ass}
\sup_{\beta \in B} \left\{ \frac{\partial \varphi}{\partial t}  + \L^\beta \varphi + f^\beta \right\} < -\eta < 0
\end{equation}
Because $\varphi \in \pb \cap C^{1,2}([0,T) \times \R^d)$ and the Hamiltonian is continuous in $(t,x)$ by (V\ref{Ass_mu_etc_continuous}), there is an open set $\mathcal{O}$ surrounding $(t_0,x_0)$ in $S_T$ where $\sup_{\beta \in B} \left\{ \frac{\partial \varphi}{\partial t}  + \L^\beta \varphi + f^\beta \right\} < -\eta/2$.

From the definition of $v^*$, there exists a sequence $(t_n,x_n) \in S_T \cap \mathcal{O}$ such that $(t_n,x_n) \to (t_0,x_0)$, $v(t_n,x_n) \to v^*(t_0,x_0)$ for $n\to\infty$. By continuity of $\varphi$, again $\delta_n := v(t_n,x_n) - \varphi(t_n,x_n)$ converges to $0$ as $n$ goes to infinity. 

By definition of the value function, there exists for all $n$ and $\varepsilon >0$ (choose $\varepsilon = \varepsilon_n$ with $\varepsilon_n \downarrow 0$) a combined admissible control $\alpha^n = \alpha^n(\varepsilon) = (\beta^n,\gamma^n)$, $\gamma^n= ( \tau^n_i, \zeta^n_i)_{i\ge 1}$ such that
\begin{equation} 
\label{subso_Jepsilon}
v(t_n, x_n) \le J^{(\alpha^n)}(t_n,x_n) + \varepsilon. 
\end{equation}

For a $\rho>0$ chosen suitably small (i.e. $B(x_n,\rho) \subset B(x_0,2\rho) \subset S$, $t_n + \rho < t_0 + 2 \rho < T$ for large $n$), we define the stopping time $\bar{\tau}_n := \bar{\tau}_n^\rho \wedge \tau_1^n$, where
\[
\bar{\tau}_n^\rho := \inf\{ s \ge t_n: | X^{\alpha^n,t_n,x_n}_s - x_n| \ge \rho \} \wedge (t_n + \rho).
\]
We want to show that $\bar{\tau}_n \to t_0$ in probability. From (\ref{subso_Jepsilon}) combined with the Markov property (\ref{dpp_ineq2}), it immediately follows that (again $\E^n = \E^{(t_n,x_n)}$, and $\beta = \beta^n$, $\alpha = \alpha^n$)\footnote{In the following, we will switch between $\alpha$ and $\beta$ in our notation, where the usage of $\beta$ indicates that there is no impulse to take into account.}
\begin{eqnarray}
v(t_n,x_n) 
&\le& \label{subso_dpp_jump} \E^n \left[ \int_{t_n}^{\bar{\tau}_n} f(s,X^\beta_s,\beta_s) ds +  v(\bar{\tau}_n,\check{X}^\beta_{\bar{\tau}_n-}) \right] + \varepsilon_n\\
&\le&  \E^n \left[ \int_{t_n}^{\bar{\tau}_n} f(s,X^\beta_s,\beta_s) ds +  \varphi(\bar{\tau}_n,\check{X}^\beta_{\bar{\tau}_n-}) \right] + \varepsilon_n.\nonumber
\end{eqnarray}
Thus by Dynkin's formula on $\varphi(\bar{\tau}_n,\check{X}^\beta_{\bar{\tau}_n-})$ and using $v(t_n,x_n) = \varphi(t_n,x_n) + \delta_n$,
\begin{eqnarray*}
\delta_n 
&\le& \E^n \left[ \int_{t_n}^{\bar{\tau}_n} f(s,X^\beta_s,\beta_s) + \frac{\partial\varphi}{\partial s}(s,X^\beta_s) + \L^\beta \varphi(s,X^\beta_s) ds \right] + \varepsilon_n\\
&\le& -\frac{\eta}{2} \E^n [ \bar{\tau}_n - t_n ] + \varepsilon_n,
\end{eqnarray*}
where for $\rho$ small enough, we have applied (\ref{subso_ass}). This implies that $\lim_{n\to\infty} \E [ \bar{\tau}_n ] = t_0$, which is equivalent to $\bar{\tau}_n \to t_0$ in probability (as one can easily check with Chebyshev's inequality; $\bar{\tau}_n$ is bounded). 

Now let us continue with our proof. In the following, we make again use of the stochastic continuity of $X^{\beta_n,t_n,x_n}$ (up until the first impulse). We define 
\[
A_n(\rho) = \{\omega: \sup_{t_n\le s\le \bar{\tau}_n} | X^{\beta_n,t_n,x_n}_s - x_n | \le \rho\}
\]
(\ref{subso_Jepsilon}) combined with the Markov property (\ref{dpp_inequality}) gives us
\begin{equation}
v(t_n,x_n) 
\le \E^n \left[ \int_{t_n}^{\bar{\tau}_n} f(s,X^\beta_s,\beta_s) ds
+ K(\tau_1^n, \check{X}^\beta_{\tau_1^n -}, \zeta_1^m ) 1_{\bar{\tau}_n^\rho \ge \tau_1^n}
+ v(\bar{\tau}_n,X^\alpha_{\bar{\tau}_n}) \right] + \varepsilon_n
\end{equation}
To find upper estimates for $v(t_n,x_n)$, we use indicator functions for three separate cases:
\begin{gather}
\tag{I} \{\bar{\tau}_n^\rho < \tau_1^n\}\\
\tag{II} \{\bar{\tau}_n^\rho \ge \tau_1^n\} \cap  A_n(\rho)^c\\
\tag{III} \{\bar{\tau}_n^\rho \ge \tau_1^n\} \cap A_n(\rho) 
\end{gather}
(III) is the predominant set: For any sequence $(\hat{\varepsilon}_n)$, by basic probability $ \P(\bar{\tau}_n^\rho \ge \tau_1^n) 
\ge 1 - \P( \bar{\tau}_n^\rho < \hat{\varepsilon}_n) - \P( \tau_1^n \ge \hat{\varepsilon}_n)$. Choose $\hat{\varepsilon}_n \downarrow t_0$ such that $\P( \tau_1^n \ge \hat{\varepsilon}_n) \to 0$. By Lemma \ref{lemma_sup_cont} (\ref{lemma_sup_cont_uniform}) (wlog, we need only consider the setting without impulses), also $\limsup_{n\to\infty} \P( \bar{\tau}_n^\rho < \hat{\varepsilon}_n) = 0$. 
In total, we have $\P(III) \to 1$ or $1_{(III)} \to 1$ a.s. for $n\to\infty$.

Thus, using that if there is an impulse in $\bar{\tau}_n$ (i.e. $\bar{\tau}_n^\rho \ge \tau_1^n$), then $v(\bar{\tau}_n,X^\alpha_{\bar{\tau}_n}) + K(\bar{\tau}_n,\check{X}^\beta_{\bar{\tau}_n-},\zeta_1^n) \le \M v(\bar{\tau}_n,\check{X}^\beta_{\bar{\tau}_n-})$,\footnote{Note: More than one impulse could occur in $\bar{\tau}_n$ if the transaction cost structure allows for it (e.g., $K$ quadratic in $\zeta$). In this case however, the result follows by monotonicity of $\M$ (Lemma \ref{lemma_Mu} (\ref{Mu_monotonous}))}
\begin{eqnarray*}
v(t_n,x_n)
&\le& \sup_{|t'-t_0| < \rho \atop |y'-x_0|< \rho} f(t',y',\beta_{t'}) \E [\bar{\tau}_n - t_n]\\
&& + \; \E [ | v(\bar{\tau}_n,X^\beta_{\bar{\tau}_n}) | 1_{(I)} ] 
+ \E [ | \M v(\bar{\tau}_n,\check{X}^\beta_{\bar{\tau}_n-})| 1_{(II)} ] \\
&& + \; \sup_{|t'-t_0| < \rho \atop |y'-x_0|< \rho} \M v(t',y') \E 1_{(III)}
\end{eqnarray*}
To prove the boundedness in term (I) (uniform in $n$), we can assume wlog that $\nu(\R^k) = 1$, and consider only jumps bounded away from $0$ for $x_0 > 0$. Then by (E\ref{ass_value_fct_in_C}), $\E \sup_n | v(\bar{\tau}_n,X^\beta_{\bar{\tau}_n}) | \le C\, \E [ 1 + R(x_0 + 2 \rho + Y) ]$ for a jump $Y$ with distribution $\nu^{\ell(t_0+ 2\rho, x_0+2 \rho,\beta,\cdot)}$, which is finite by the definition of $\pb$. The same is true for $\M v(\bar{\tau}_n,\check{X}^\beta_{\bar{\tau}_n-})$ (for $\M v \ge 0$ because $\M v \le v$, for negative $\M v$ this follows from the definition).

Sending $n\to\infty$ ($\limsup_{n\to\infty}$), the $f$-term, and term (I) converge to $0$ by the dominated convergence theorem. For term (II), a general version of the DCT shows that it is bounded by 
\[
\E [\limsup_{n\to\infty} |\M v(\bar{\tau}_n,\check{X}^\beta_{\bar{\tau}_n-})| 1_{A_n(\rho)^c} ],
\]
and term (III) becomes $\sup \M v(t',y')$.
Now we let $\rho\to 0$ and obtain:
\[
v^*(t_0,x_0) \le \lim_{\rho \downarrow 0} \sup_{|t'-t_0| < \rho \atop |y'-x_0|< \rho} \M v(t',y')
= (\M v)^*(t_0,x_0) \le \M v^*(t_0,x_0)
\]
by Lemma \ref{lemma_Mu} (\ref{Mu_usc}), a contradiction. Thus (\ref{subso_0}) is true. 
\hfill$\Box$\\

Let us elaborate on some details of the proof:
\begin{itemize}
\item In the proof, we have only used that all constant controls with values in $B$ are admissible for the SDE (\ref{viso_SDE}). So actually, we are quite free how to choose the set of admissible controls -- the value function always turns out to be a viscosity solution.
\item Another approach for the subsolution part would be tempting, although we do not see how this can work: In the subsolution proof, we assumed $v^*(t_0,x_0)- \M v^*(t_0,x_0)>0$. This implies, using Lemma \ref{lemma_Mu} (\ref{Mu_star_jump_m}) and the $0$-$1$ law, that for $n$ large enough, $\tau_1^n>t_n$ a.s. On the other hand, from $\bar{\tau}_n \to t_0$, it follows by Lemma \ref{lemma_taumrho_not0} that $\tau_1^n \to t_0$ in probability (it is sufficient to consider the setting of Lemma \ref{lemma_taumrho_not0} without impulse, since otherwise the first impulse would anyhow converge to $0$ in probability). So the convergence of $\tau_1^n$ points already to a contradiction.
\item By local boundedness of $v$, the derivation of the dynamic programming principle and the integrability condition (\ref{eq_integrability}), we could already deduce that $\E^{t_n,x_n} [ v(\bar{\tau}_n,X^\beta_{\bar{\tau}_n}) ] < \infty$; it is however not so easy to deduce this uniformly in $n$. The condition (\ref{eq_integrability}) contains implicitly conditions on $\nu$ we have formalised in (V\ref{ass_value_fct_in_C}). 
\end{itemize}
\smallskip

We promised to come back to the ``regularity at $\partial S$'' issue, and present here conditions sufficient for condition (E\ref{ass_ex_v_cont_on_boundary}).

\begin{enumerate}[(E1*)]
\item \label{ass_ex_regular_on_boundary} For any point $(t,x) \in [0,T) \times \partial S$, any sequence $(t_n,x_n) \subset [0,T) \times S$, $(t_n,x_n) \to (t,x)$, and for each small $\varepsilon>0$, there is an admissible combined control $\alpha_n = (\beta_n,\gamma_n)$ such that
\begin{equation}
v(t_n,x_n) \le J^{(\alpha_n)}(t_n,x_n) + \varepsilon,
\end{equation}
and such that for all $\delta>0$, $\P(\tilde{\tau}_S^n < t_n + \delta) \to 1$ for $n \to \infty$ (where $\tilde{\tau}_S^n = \inf\{s\ge t_n: X^{\beta_n,t_n,x_n}(s) \not\in S \}$).
\item \label{ass_ex_NIR_open_on_boundary} For any point $(t,x) \in \partial^*S_T$, if there is a sequence $(t_n,x_n) \subset [0,T) \times S$ converging to $(t,x)$ with $v(t_n,x_n) > \M v(t_n,x_n)$, then there is a neighbourhood of $(t,x)$ (open in $[0,T] \times \R^d$), where $v > \M v$
\end{enumerate}

\begin{example}
Let $X$ be a one-dimensional Brownian motion with $\sigma>0$, and assume it is never optimal to give an impulse near the boundary. Then (E\ref{ass_ex_regular_on_boundary}*) and (E\ref{ass_ex_NIR_open_on_boundary}*) are satisfied.
\end{example}

We define $\tau_S^n = \inf\{s\ge t_n: X^{\alpha_n,t_n,x_n}(s) \not\in S \}$.

\begin{proposition}
Let $\tau = \tau_S$ or $\tau = \tau_S \wedge T$. If (E\ref{ass_ex_regular_on_boundary}*) and (E\ref{ass_ex_NIR_open_on_boundary}*) hold, and for $n$ large the integrability condition 
\[
\int_t^{\tau} |f(s,X^{\alpha_n}_s,(\beta_n)_s)| ds + |g(\tau, X^{\alpha_n}_{\tau})| 1_{\tau < \infty} + \sum_{\tau_j \le \tau} |K(\tau^n_j, \check{X}^{\alpha_n}_{\tau^n_j -}, \zeta^n_j )| \le h \in L^1(\P;\R)
\]
is satisfied, then (E\ref{ass_ex_v_cont_on_boundary}) holds.
\end{proposition}

\textbf{Proof: } Let $(t,x) \in [0,T) \times \partial S$ and assume wlog that $v > \M v$ in a neighbourhood of $(t,x)$. Then we have to show for all $\varepsilon >0$ and all sequences chosen as in (E\ref{ass_ex_regular_on_boundary}*), that
\begin{equation}
\label{eq_regular_boundary_convergence}
\E^{(t_n,x_n)} \left[ \int_{t_n}^{\tau_S} f(s,X^{\alpha_n}_s,(\beta_n)_s) ds + g(\tau_S,X^{\alpha_n}(\tau_S)) 1_{\tau_S < \infty} + \sum_{\tau_j \le \tau_S} K(\tau^n_j, \check{X}^{\alpha_n}_{\tau^n_j -}, \zeta^n_j ) \right] \to g(t,x)
\end{equation}
as $n \to \infty$, $\delta \to 0$. For the set 
\[
B_{n,\delta} := \{ \tau_S^n < t_n + \delta \} \cap \{ \sup_{t_n \le t \le \tau_S^n} | X_t^{\alpha_n,t_n,x_n} - x_n | < \delta \},
\]
we claim that for all $\delta>0$ small enough, $1_{B_{n,\delta}} \to 1$ as $n\to\infty$. To see this, first note that by assumption and Lemma \ref{lemma_sup_cont} (\ref{lemma_sup_cont_uniform}), this is true for the set
\[
\tilde{B}_{n,\delta} := \{ \tilde{\tau}_S^n < t_n + \delta \} \cap \{ \sup_{t_n \le t \le \tilde{\tau}_S^n} | X_t^{\beta_n,t_n,x_n} - x_n | < \delta \}.
\]
Choose $\delta$ small enough such that $v > \M v$ on $B(x_n,\delta)$ and $x \in B(x_n,\delta)$ for all $n$ large. Then it is easily checked that $\{ \sup_{t_n \le t \le \tau_S^n} | X_t^{\alpha_n,t_n,x_n} - x_n | < \delta \} = \{ \sup_{t_n \le t \le \tilde{\tau}_S^n} | X_t^{\beta_n,t_n,x_n} - x_n | < \delta \}$. 
Because 
\[
\P\left( \tau_S^n < t_n + \delta \left| \sup_{t_n \le t \le \tilde{\tau}_S^n} | X_t^{\beta_n,t_n,x_n} - x_n | < \delta \right. \right) 
= \P\left( \tilde{\tau}_S^n < t_n + \delta \left| \sup_{t_n \le t \le \tilde{\tau}_S^n} | X_t^{\beta_n,t_n,x_n} - x_n | < \delta \right. \right) \to 1,
\]
we can conclude that $1_{B_{n,\delta}} \to 1$ as $n\to\infty$ for $\delta>0$ small enough.

The convergence in (\ref{eq_regular_boundary_convergence}) then follows just as in the existence proof, by splitting into $B_{n,\delta}$ and $B_{n,\delta}^c$ and applying the $\limsup$, $\liminf$ versions of the DCT, first for $n\to \infty$, and then for $\delta \to 0$ (using the continuity of $f,g$). The result for $t=T$ holds by the same arguments, because time is always regular.
\hfill$\Box$

\begin{rem}
(E\ref{ass_ex_v_cont_on_boundary}) (resp, (E\ref{ass_ex_regular_on_boundary}*), (E\ref{ass_ex_NIR_open_on_boundary}*)) excludes in particular problems with \emph{de facto} state constraints, where it is optimal to stay inside $S$. We note however that the framework presented here allows for an adaptation to (true and \emph{de facto}) state constraints, which can be pretty straightforward for easy constraints. Apart from the stochastic proof that we can restrain ourselves to controls keeping the process inside $S$, the adaptation involves changing the function $w$ used in the uniqueness part, such that only values in $S$ need to be considered in the comparison proof. For an example in the diffusion case, see \citet{LVMP07}; jumps outside $S$ however may be difficult to handle. 
\end{rem}


\subsection{Elliptic case}
\label{subsec_existence_elliptic}

The existence result for the elliptic QVI (\ref{qvi_elliptic}) now follows from the parabolic result by an exponential time transformation. Recall that the elliptic QVI is
\begin{equation}
\tag{\ref{qvi_elliptic}}
\min ( - \sup_{\beta \in B} \{ \L^\beta u + f^\beta \}, u - \M u ) = 0,
\end{equation}
for the elliptic integro-differential operator $\L^\beta$ from (\ref{L_beta_elliptic}) and the intervention operator $\M$ selecting the optimal instantaneous impulse. 

For an $s\ge 0$ and $\rho>0$, let the functions as used in section \ref{subsec_existence_parabolic} be tagged by a tilde. $\tilde{f}$, $\tilde{g}$, $\tilde{K}$ are all built in the same way on $[0,\infty) \times \R^d$, as the following example: 
\[
\tilde{f}(t,x,\beta) = e^{-\rho(s+t)} f(x,\beta)
\]
Let $\tilde{\Gamma} = \Gamma$ (independent of $t$) and likewise $\tilde{\mu},\tilde{\sigma},\tilde{\ell}$ and the transaction set $\tilde{Z}$. The intervention operator including time (as in section \ref{subsec_existence_parabolic}) is denoted by $\tilde{\M}$, and the time-independent one is defined by
\[
\M u(x) = \sup \{ u(\Gamma(x,\zeta)) + K(x,\zeta): \zeta \in Z(x) \},
\]
Let $\tilde{\L}^\beta$ denote the integro-differential operator (or infinitesimal generator) from (\ref{L_beta_parabolic}) as used in section \ref{subsec_existence_parabolic} (in this case $\tilde{\L}^\beta$ is not equal to $e^{-\rho(s+t)} \L^\beta$). 

It can be checked that the assumptions of section \ref{subsec_existence_parabolic} hold for the tilde functions, if the corresponding assumption holds for the time-independent functions without tilde. As well, all the lemmas used for the proof of the existence theorem are still valid in the time-independent case.

\begin{definition}[Viscosity solution]
A function $u \in \pb(\R^d)$ is a (viscosity) subsolution of (\ref{qvi_elliptic}) if for all $x_0 \in \R^d$ and $\varphi \in \pb \cap C^2(\R^d)$ with $\varphi(x_0) = u^*(x_0)$, $\varphi \ge u^*$,
\begin{eqnarray*}
\min \left( - \sup_{\beta \in B} \left\{ \L^\beta \varphi + f^\beta \right\}, u^* - \mathcal{M}u^* \right) 
&\le & 0 \qquad \mbox{ in } x_0 \in S\\
\min \left( u^* - g, u^* - \mathcal{M}u^* \right) 
&\le & 0 \qquad \mbox{ in } x_0 \in \R^d \setminus S
\end{eqnarray*}
A function $u \in \pb(\R^d)$ is a (viscosity) supersolution of (\ref{qvi_elliptic}) if for all $x_0 \in \R^d$ and $\varphi \in \pb \cap C^2(\R^d)$ with $\varphi(x_0) = u_*(x_0)$, $\varphi \le u_*$,
\begin{eqnarray*}
\min \left( - \sup_{\beta \in B} \left\{ \L^\beta \varphi + f^\beta \right\}, u_* - \mathcal{M}u_* \right) 
&\ge & 0 \qquad \mbox{ in } x_0 \in S\\
\min \left( u_* - g, u_* - \mathcal{M}u_* \right) 
&\ge & 0 \qquad \mbox{ in } x_0 \in \R^d \setminus S
\end{eqnarray*}
A function $u$ is a viscosity solution if it is sub- and supersolution.
\end{definition}

In the original problem on $[0,\infty) \times \R^d$, we only consider Markov controls that are time-independent, i.e., only depend on the state variable $x$. Denote by $\tilde{v}$ the then resulting value function of the parabolic problem on $[0,\infty) \times \R^d$. Then we can define the time-independent value function 
\begin{equation}
v(x) := e^{\rho(s+t)} \tilde{v}(t,x).
\end{equation}
Let us emphasise that this definition is only admissible if the right-hand side actually does not depend on $t$, which can be checked in the definition of $J^{(\alpha)}(t,x)$ by the homogeneous Markov property (and of course only if the time horizon is infinite, i.e., $\tau = \tau_S$ for $S \subset \R^d$).

The existence of the elliptic QVI then follows by an easy time transformation:

\begin{corollary}
\label{theo_existence_qvi_elliptic}
Let Assumptions \ref{ass_existence_uniqueness} and \ref{ass_theo_existence_qvi} be satisfied. Then the value function $v$ as defined above is a viscosity solution of (\ref{qvi_elliptic}).
\end{corollary}

\textbf{Proof: } We know from Theorem \ref{theo_existence_qvi_parabolic} that $\tilde{v}$ is a (parabolic) viscosity solution of (\ref{qvi_parabolic}) on $[0,\infty) \times \R^d$ (without the terminal condition), i.e. of
\begin{equation}
\begin{split}
\min ( - \sup_{\beta \in B} \{ u_t + \tilde{\L}^\beta u + \tilde{f}^\beta \}, u - \tilde{\M}u ) &= 0 \qquad \mbox{ in } (t_0,x_0) \in [0,\infty) \times S\\
\min( u - \tilde{g}, u - \tilde{\M} u ) &= 0 \qquad \mbox{ in } (t_0,x_0) \in [0,\infty) \times (\R^d \setminus S)
\end{split}
\end{equation}
For the subsolution proof, let $\varphi \in \pb \cap C^2(\R^d)$ with $\varphi(x_0) = v^*(x_0)$, $\varphi \ge v^*$. For $\tilde{\varphi}(t,x) := e^{-\rho(s+t)} \varphi(x)$ for all $t\ge 0$, $\tilde{\varphi}(t_0,x_0) = \tilde{v}^*(t_0,x_0)$ and $\tilde{\varphi} \ge \tilde{v}^*$. Furthermore, 
\[
\tilde{\varphi}_t + \tilde{\L}^\beta \tilde{\varphi} 
= e^{-\rho(s+t)} ( -\rho \varphi + e^{\rho(s+t)} \tilde{\L}^\beta \tilde{\varphi} )
= e^{-\rho(s+t)} \, \L^\beta \varphi 
\]
By the definition of the elliptic $\M$, we have $\tilde{v} - \tilde{\M}\tilde{v} = e^{-\rho(s+t)} ( v - \M v)$. The supersolution property is proved in the same manner.
\hfill$\Box$


\section{Uniqueness}
\label{sec_uniqueness}

The purpose of this section is to prove uniqueness results both for the elliptic and the parabolic QVI by analytic means. Using such a uniqueness result, together with the existence results of section \ref{sec_existence}, we can conclude
\begin{quote}
\begin{center} The viscosity solution of the QVI is equal to the value function. \end{center}
\end{quote}
We were inspired mainly by the papers \citet{Is93} (for the impulse part) and \citet{BaIm07} (for the PIDE part). As general reference for viscosity solutions, \citet{CIL} was used and will be frequently cited. Some ideas have come from \cite{OkSu02}, \cite{OkSu05}, \cite{AkMeSu96} and \cite{JaKa06}.

In this section, $v$ does not denote the value function any longer, and some other variables may serve new purposes as well.

First, we will investigate uniqueness of QVI viscosity solutions for the elliptic case (the parabolic case following at the end):
\begin{equation}
\tag{\ref{qvi_elliptic}}
\begin{split}
\min ( - \sup_{\beta \in B} \{ \L^\beta u + f^\beta \}, u - \M u ) &= 0 \qquad \mbox{in } S\\
\min( u - g, u - \M u ) &= 0 \qquad \mbox{in } \R^d \setminus S
\end{split}
\end{equation}


\subsection{Preliminaries}
\label{subsec_uniqueness_prelims}

Whereas in the last section, we did not care about the specific form of the generator (as long as the Dynkin formula was valid), we now need to investigate the operator $\L^\beta$ more in detail:
\begin{multline}
\label{L_beta_detailed}
\L^\beta u (x) = \frac{1}{2} tr \left(\sigma(x,\beta)\sigma^T(x,\beta)D^2 u(x)\right) + \langle \mu(x,\beta), \grad u(x) \rangle - c(x) u(x) \\
+ \int u(x + \ell(x,\beta,z)) - u(x) - \langle \grad u(x), \ell(x,\beta,z) \rangle 1_{|z|< 1} \, \nu(dz),
\end{multline}
where $c$ is some positive function related to the discounting in the original model.\footnote{Note that the normal definition of a Lévy integral is with the indicator function $1_{|z|\le 1}$; this is however equivalent}

We recall the definition of the function space $\pb = \pb(\R^d)$ from section \ref{sec_model}, such that the differential operator $L^\beta$ is well-defined for $\phi \in \pb \cap C^2(\R^d)$.  
Denoting for $0< \delta < 1$, $y,p\in \R^d$, $X\in \R^{d \times d}$, $r \in \R$ and $(l_\beta)_{\beta \in B} \subset \R$:
\begin{align*}
F(x,r,p,X,(l_\beta)) &= - \sup_{\beta \in B} \left\{ \frac{1}{2} tr \left(\sigma(x,\beta)\sigma^T(x,\beta) X \right) + \langle \mu(x,\beta), p \rangle - c(x) r + f(x,\beta) + l_\beta \right\}\\
\I^{1,\delta}_\beta [x,\phi] &= \int_{|z|<\delta} \phi(x + \ell(x,\beta,z)) - \phi(x) - \langle \grad \phi(x), \ell(x,\beta,z) \rangle  \, \nu(dz)\\
\I^{2,\delta}_\beta [x,p,\phi] &= \int_{|z|\ge \delta} \phi(x + \ell(x,\beta,z)) - \phi(x) - \langle p, \ell(x,\beta,z) \rangle 1_{|z|< 1}  \, \nu(dz)\\
\I_\beta[x,\phi] &=  \I^{1,\delta}_\beta [x,\phi] + \I^{2,\delta}_\beta [x,\grad \phi(x),\phi],
\end{align*}
we have to analyse the problem 
\[
\min ( F(x,u(x),\grad u(x),D^2 u(x), \I_\beta [x, u(\cdot)] ), u(x) - \M u(x) ) = 0,
\]
where the notation $u(\cdot)$ in the integral indicates that nonlocal terms are used on $u$, not only from $x$. As well, $\I_\beta [x, u(\cdot)]$ within $F$ always stands for a family ($\beta \in B$) of integrals. Denote by $F^\beta$ the function $F$ without the sup, i.e. for a concrete $\beta$.

\begin{rem}
The following properties hold for our problem:
\begin{enumerate}[(P1)]
\item \label{Prop_ellipticity} Ellipticity of $F$: $F(x,r,p,X^1,(l^1_\beta)) \le F(x,r,p,X^2,(l^2_\beta))$ if $X^1 \ge X^2$, $l^1_\beta \ge l^2_\beta$ $\forall \beta \in B$
\item Translation invariance: $u - \M u = (u+l) - \M (u+l)$, $\I[y_0,\phi] = \I[y_0,\phi + l]$ for constants $l \in \R$
\item \label{Prop_continuity_l} $(l_\beta)_\beta \mapsto F(x,r,p,X,(l_\beta))$ is continuous in the sense that
\[
|F(x,r,p,X,(l^1_\beta)) - F(x,r,p,X,(l^2_\beta)) | \le \sup_\beta | l^1_\beta - l^2_\beta |.
\]
\end{enumerate}
\end{rem}
The last statement -- proved by easy sup manipulations -- is just for the sake of completeness; we will not use it explicitly because uniform convergence needs continuous functions, which in general we do not have. 

For $x \in S$ and $\delta>0$, recall the definition of $U_\delta = U_\delta(x) = \{ \ell(x,\beta,z): |z| < \delta \}$. $U_\delta$ facilitates splitting the integral $\I_\beta[x,\phi]$: changing $\phi$ on $U_\delta(x)$ only affects $\I^{1,\delta}_\beta [x,\phi]$, and reversely.

Let henceforth Assumptions \ref{ass_existence_uniqueness} and \ref{ass_uniqueness_qvi} be satisfied (the latter needed mainly for the equivalence of different viscosity solution definitions). Further assume 
\begin{enumerate}[(U1*)]
\item \label{U_ass_c_cont} $c$ is continuous 
\end{enumerate}

\begin{rem}
It is sufficient for the comparison theorem if Assumption \ref{ass_uniqueness_qvi} holds only for small $\delta>0$: ``\textit{For any $x_0$, there is a small environment and a $\bar{\delta}>0$, where the assumption holds for $0 < \delta < \bar{\delta}$\ldots}.'' This is why we will carry out all proofs for $\ell$ depending on $\beta$ in the following.
\end{rem}

Immediately from (V\ref{Ass_mu_etc_continuous}) and (U\ref{U_ass_c_cont}*), it follows that $\sup_{\beta \in B} \sigma(x,\beta) < \infty$, and by sup manipulations that $(x,r,p,X) \mapsto F(x,r,p,X,l_\beta)$ is continuous; but even more can be deduced:

\begin{proposition}
\label{prop_uniqueness_F_continuity}
Let $(\beta_k) \subset B$ with $\beta_k \to \beta$, and $(x_k),(p_k) \subset \R^d$ with $x_k \to x \in S$, $p_k \to p$. 
\begin{enumerate}[(i)]
\item \label{prop_I2_usc} If $u \in \pb \cap USC(\R^d)$ and $v \in \pb \cap LSC(\R^d)$ with $u(x_k) \to u(x)$, $v(x_k) \to v(x)$, then
\begin{align*}
\limsup_{k\to\infty} \I^{2,\delta}_{\beta_k}[x_k,p_k,u(\cdot)] &\le \I^{2,\delta}_\beta[x,p,u(\cdot)],
& \liminf_{k\to\infty} \I^{2,\delta}_{\beta_k}[x_k,p_k,v(\cdot)] &\ge \I^{2,\delta}_\beta[x,p,v(\cdot)].
\end{align*}
Moreover, for $(\varphi_k), (\psi_k) \subset \pb \cap C^2(\R^d)$ with $\varphi_k \to u$ and $\psi_k \to v$ monotonously, $\varphi_k(x_k) \to u(x)$, $\psi_k(x_k) \to v(x)$,
\begin{align*}
\limsup_{k\to\infty} \I^{2,\delta}_{\beta_k}[x_k,p_k,\varphi_k(\cdot)] &\le \I^{2,\delta}_\beta[x,p,u(\cdot)],
& \liminf_{k\to\infty} \I^{2,\delta}_{\beta_k}[x_k,p_k,\psi_k(\cdot)] &\ge \I^{2,\delta}_\beta[x,p,v(\cdot)].
\end{align*}
\item \label{prop_I1_cont} If $\varphi \in C^2(\R^d)$, then $(x,\beta) \mapsto \I^{1,\delta}_\beta[x,\varphi(\cdot)]$ is continuous. Moreover, for $(\varphi_k) \subset C^2(\R^d)$ with $\varphi_k \to \varphi$ monotonously, $\varphi_k = \varphi$ in an environment of $x$,
\[
\lim_{k\to\infty} \I^{1,\delta}_{\beta_k}[x_k,\varphi_k(\cdot)] = \I^{1,\delta}_\beta[x,\varphi(\cdot)].
\]
\item \label{prop_F_lsc_in_phi} If $u \in \pb \cap USC(\R^d)$ and $v \in \pb \cap LSC(\R^d)$, $(\varphi_k), (\psi_k) \subset \pb \cap C^2(\R^d)$ with $\varphi_k \to u$ and $\psi_k \to v$ monotonously, $\varphi_k = u \in C^2$ and $\psi_k = v \in C^2$ in an environment of $x$, then
\begin{align}
\label{eq_prop_F_lsc_in_phi_liminf}
\liminf_{k\to\infty} F^{\beta_k}(x,r,p,X,\I_{\beta_k}[x_k,\varphi_k(\cdot)]) 
&\ge F^{\beta}(x,r,p,X,\I_{\beta}[x,u(\cdot)])\\
\label{eq_prop_F_lsc_in_phi_limsup}
\limsup_{k\to\infty} F^{\beta_k}(x,r,p,X,\I_{\beta_k}[x_k,\psi_k(\cdot)]) 
&\le F^{\beta}(x,r,p,X,\I_{\beta}[x,v(\cdot)])
\end{align}
\item \label{prop_F_usc_in_beta} If $u \in \pb \cap USC(\R^d)$, $\varphi \in C^2(\R^d)$, then $\beta \mapsto -F^\beta(x,r,p,X,\I^{1,\delta}_\beta[x,\varphi(\cdot)] + \I^{2,\delta}_\beta[x,p,u(\cdot)]$ is in $USC(\R^d)$. In particular, the supremum in $\beta \in B$ is assumed.
\end{enumerate}
\end{proposition}

\textbf{Proof: } (\ref{prop_I2_usc}): We prove only the first statement for $u \in USC$ (the $lsc$ proof being analogous). By a general version of the DCT and the definition of $\pb$,  
\begin{eqnarray}
&& \limsup_{k\to\infty} \int_{|z|\ge \delta} u(x_k + \ell(x_k,\beta_k,z)) - u(x_k) - \langle p_k, \ell(x_k,\beta_k,z) \rangle 1_{|z|< 1}  \, \nu(dz)\nonumber\\
&\le& \int_{|z|\ge \delta} \limsup_{k\to\infty} \, u(x_k + \ell(x_k,\beta_k,z)) - u(x_k) - \langle p_k, \ell(x_k,\beta_k,z) \rangle 1_{|z|\le 1}  \, \nu(dz)\nonumber\\
\label{cor_liminf_limsup_int} 
&\le& \int_{|z|\ge \delta} u(x + \ell(x,\beta,z)) - u(x) - \langle p, \ell(x,\beta,z) \rangle 1_{|z|< 1}  \, \nu(dz),
\end{eqnarray}
where we have used the continuity of $\ell$ (V\ref{Ass_mu_etc_continuous}). The fact for $\varphi_k$ follows because by monotonicity of $(\varphi_k)$, and a general version of Dini's Theorem (cf. \citet{Di02}, Th. 7.3), $\limsup_{k \to \infty} \varphi_k(x_k + \ell(x_k,\beta_k,z)) \le u^*(x + \ell(x,\beta,z))$. 
\\
(\ref{prop_I1_cont}): Outside of the singularity of $\nu$, the result follows from (\ref{prop_I2_usc}). In the environment of $x$, the Taylor expansion (\ref{I1_Taylor_expansion}) gives the upper bound for the application of the DCT (where the local boundedness of $U_\delta(x)$ holds by (V\ref{Ass_ell_integrable}).\\
(\ref{prop_F_lsc_in_phi}) follows immediately from (\ref{prop_I2_usc}), (\ref{prop_I1_cont}). Finally, (\ref{prop_F_usc_in_beta}) holds by the continuity conditions (V\ref{Ass_mu_etc_continuous}), (U\ref{U_ass_c_cont}*) and (\ref{prop_I2_usc}), (\ref{prop_I1_cont}). \hfill$\Box$


\subsection{Viscosity solutions: Different definitions}
\label{subsec_uniqueness_defs}

Let us now restate in the new notation the definition of viscosity solution in the elliptic case (equivalent thanks to the translation invariance property):

\begin{definition}[Viscosity solution 1]
\label{viso_def1}
A function $u \in \pb$ is a (viscosity) subsolution of (\ref{qvi_elliptic}) if for all $x_0 \in \R^d$ and $\varphi \in \pb \cap C^2(\R^d)$ such that $u^* - \varphi$ has a global maximum in $x_0$,
\begin{eqnarray*}
\min \left( F(x_0,u^*,\grad \varphi,D^2 \varphi, \I_\beta [x_0, \varphi(\cdot)] )
, u^* - \mathcal{M}u^* \right) 
&\le & 0 \qquad \mbox{ in } x_0 \in S\\
\min \left( u^* - g, u^* - \mathcal{M}u^* \right) 
&\le & 0 \qquad \mbox{ in } x_0 \in (\R^d \setminus S)
\end{eqnarray*}
A function $u \in \pb$ is a (viscosity) supersolution of (\ref{qvi_elliptic}) if for all $x_0 \in \R^d$ and $\varphi \in \pb \cap C^2(\R^d)$ such that $u_* - \varphi$ has a global minimum in $x_0$,
\begin{eqnarray*}
\min \left( F(x_0,u_*,\grad \varphi,D^2 \varphi, \I_\beta [x_0, \varphi(\cdot)] )
, u_* - \mathcal{M}u_* \right) 
&\ge & 0 \qquad \mbox{ in } x_0 \in S\\
\min \left( u_* - g, u_* - \mathcal{M}u_* \right) 
&\ge & 0 \qquad \mbox{ in } x_0 \in (\R^d \setminus S)
\end{eqnarray*}
A function $u$ is a viscosity solution if it is sub- and supersolution.
\end{definition}

As in \cite{BaIm07}, we give two equivalent definitions, the latter one is needed later on in the uniqueness proof:

\begin{definition}[Viscosity solution 2]
\label{viso_def2}
A function $u \in \pb$ is a (viscosity) subsolution of (\ref{qvi_elliptic}) if for all $x_0 \in \R^d$ and $\varphi \in C^2(\R^d)$ such that $u^* - \varphi$ has a maximum in $x_0$ on $U_\delta(x_0)$,
\begin{eqnarray*}
\min \left( F(x_0,u^*,\grad \varphi,D^2 \varphi, \I_\beta^{1,\delta} [x_0, \varphi(\cdot)] + \I_\beta^{2,\delta} [x_0,\grad \varphi, u^*(\cdot)]  )
, u^* - \mathcal{M}u^* \right) 
&\le & 0 \qquad \mbox{ in } x_0 \in S\\
\min \left( u^* - g, u^* - \mathcal{M}u^* \right) 
&\le & 0 \qquad \mbox{ in } x_0 \in (\R^d \setminus S)
\end{eqnarray*}
A function $u \in \pb$ is a (viscosity) supersolution of (\ref{qvi_elliptic}) if for all $x_0 \in \R^d$ and $\varphi \in C^2(\R^d)$ such that $u_* - \varphi$ has a minimum in $x_0$ on $U_\delta(x_0)$,
\begin{eqnarray*}
\min \left( F(x_0,u_*,\grad \varphi,D^2 \varphi, \I_\beta^{1,\delta} [x_0, \varphi(\cdot)] + \I_\beta^{2,\delta} [x_0,\grad \varphi, u_*(\cdot)] )
, u_* - \mathcal{M}u_* \right) 
&\ge & 0 \qquad \mbox{ in } x_0 \in S\\
\min \left( u_* - g, u_* - \mathcal{M}u_* \right) 
&\ge & 0 \qquad \mbox{ in } x_0 \in (\R^d \setminus S)
\end{eqnarray*}
A function $u$ is a viscosity solution if it is sub- and supersolution.
\end{definition}

Note that the definition is of course still valid if $U_\delta(x_0) \cap (\R^d \setminus S) \not= \emptyset$. That the two definitions are equivalent is essentially not a new result; for proofs in simpler settings, see \citet{AlTo96}, \citet{JaKa06}. 

\begin{proposition}
\label{prop_eq_def12}
Definitions 1 and 2 are equivalent.
\end{proposition}

\textbf{Proof: } 
We treat first the subsolution part. Because $u^* - \mathcal{M}u^* \le 0$ holds independently of the viscosity formulation, so of the case distinction, we need only consider the PIDE part for $x_0 \in S$. (Note that for $\nu(\R^k)<\infty$, we can set $\delta=0$, and use in the following an arbitrary local environment of $x_0$ instead of $U_\delta(x_0)$.)
\begin{itemize}
\item[``$\Leftarrow$''] Assume $u$ is a viscosity subsolution according to Definition 2. Let $x_0 \in S$ and $\varphi \in \pb \cap C^2(\R^d)$ such that $u^* - \varphi$ has a global maximum in $x_0$. Then $x_0$ is also a maximum point on $U_\delta(x_0)$. So 
\[
F( x_0,u^*,\grad \varphi,D^2 \varphi, \I_\beta^{1,\delta} [x_0, \varphi(\cdot)] + \I_\beta^{2,\delta} [x_0,\grad \varphi, u^*(\cdot)] ) \le 0
\]
By the property $\I_\beta^{2,\delta} [x_0,\grad \varphi, u^*(\cdot)] \le \I_\beta^{2,\delta} [x_0,\grad \varphi, \varphi(\cdot)]$, and ellipticity of $F$, we are done. 
 
\item[``$\Rightarrow$''] Assume $u$ is a viscosity subsolution according to Definition 1. Let $x_0 \in S$ and $\varphi \in C^2(\R^d)$ such that $u^* - \varphi$ has a maximum in $x_0$ on $U_\delta(x_0)$. (Wlog assume that $u^*(x_0) = \varphi(x_0)$, $u^* \le \varphi$ on $U_\delta(x_0)$.) 

Now consider the function
\[
\psi(x) := 1_{\R^d \setminus U_\delta}(x) u^*(x) + 1_{U_\delta}(x) \varphi(x),
\]
which is in $\pb$. It is immediate that $\psi$ is upper semicontinuous if $U_\delta$ is closed, and in this case, we can construct a monotonously decreasing sequence (e.g., by approximating with piecewise constant functions, smoothed by the standard mollifier) 
$(\varphi_k) \subset \pb \cap C^2(\R^d)$ such that $\varphi_k \downarrow \psi$ pointwise. 

If $U_\delta$ is not closed, then $\psi$ need not be upper semicontinuous (if $u^*(y) < \varphi(y)$ for a $y \in \partial U_\delta \cap (U_\delta)^c$). In this case however, there is a small open neighourhood $V$ of $y$ where $u^* < \varphi$ (because $\varphi \in LSC$, $u^* \in USC$). So we can approximate $\varphi$ from below in $V \cap int(U_\delta)$ with $\varphi_k \ge u^*$, and combine this in the construction as above. 

Let further $\varphi_k = \varphi$ in $B(x_0,\rho) \subset\subset U_\delta(x_0)$ for some $\rho>0$ (possible by (U\ref{U_ass_Upositive})). Then $u^*(x_0)=\varphi_k(x_0)$, $u^* \le \varphi_k$, and all local properties for $\varphi_k$ in $x_0$ are inherited from $\varphi$. Applying Definition 1 yields
\[
F( x_0,u^*,\grad \varphi,D^2 \varphi, \I_\beta^{1,\delta} [x_0, \varphi_k(\cdot)] + \I_\beta^{2,\delta} [x_0,\grad \varphi, \varphi_k(\cdot)] ) \le 0
\]
For each $k$, by Proposition \ref{prop_uniqueness_F_continuity} (\ref{prop_F_usc_in_beta}), the supremum in $F$ is attained in a $\beta_k$. Choose a subsequence such that $\beta_k \to \beta$. Using the limit in (\ref{eq_prop_F_lsc_in_phi_liminf}), we obtain our result. 

\end{itemize}

Supersolution part, ``$\Rightarrow$'': For analogously defined $\varphi$, here we have to show that 
\[
F( x_0,u_*,\grad \varphi,D^2 \varphi, \I_\beta^{1,\delta} [x_0, \varphi(\cdot)] + \I_\beta^{2,\delta} [x_0,\grad \varphi, u_*(\cdot)] ) \ge 0
\]
for which it is sufficient to prove $F^\beta(\ldots) \ge 0$ for all $\beta \in B$. We construct the sequence $(\varphi_k)$ with $\varphi_k \to \psi$,
\[
\psi(x) := 1_{\R^d \setminus U_\delta}(x) u_*(x) + 1_{U_\delta}(x) \varphi(x),
\]
as in the subsolution case (where the problematic part is $u_*(y) > \varphi(y)$ in a $y \in \partial U_\delta \cap (U_\delta)^c$). The result is obtained analogously (with $\beta$ fixed) to the subsolution case using (\ref{eq_prop_F_lsc_in_phi_limsup}).
\hfill$\Box$

\begin{rem}
The approximation from above of a usc function by $C^2$ functions used in the second part of the proof of Proposition \ref{prop_eq_def12} shows also that, in Definitions 1 and 2, we may restrict ourselves to functions $\varphi \in C^0(\R^d)$ that are only $C^2$ in a small neighbourhood of $x_0$.
\end{rem}

We recall the semijets needed for a third equivalent definition. They are motivated by a classical property of differentiable functions. Let $u:\R^d \to \R$.
\begin{align*}
J^+u(x) &= \{ (p,X) \in \R^d \times \spd^d: u(x+z) \le u(x) + \langle p, z \rangle + \frac{1}{2} \langle X z, z \rangle + o(|z|^2) \mbox{ as } z \to 0 \}\\
J^-u(x) &= \{ (p,X) \in \R^d \times \spd^d: u(x+z) \ge u(x) + \langle p, z \rangle + \frac{1}{2} \langle X z, z \rangle + o(|z|^2) \mbox{ as } z \to 0 \}
\end{align*}
If $u$ is twice differentiable at $x$, then $J^+u(x) \cap J^-u(x) = \{(\grad u(x), D^2 u(x)) \}$. The limiting semijets are defined by, e.g.,
\begin{multline*}
\overline{J}^+u(x) = \{ (p,X) \in \R^d \times \spd^d: \mbox{ there exist } (x_k,p_k,X_k) \to (x,p,X), \\
(p_k,X_k) \in J^+u(x_k) \mbox{ 	such that } u(x_k) \to u(x) \}.
\end{multline*}

\begin{definition}[Viscosity solution 3]
\label{viso_def3}
A function $u \in \pb$ is a (viscosity) subsolution of (\ref{qvi_elliptic}) if for all $x_0 \in \R^d$ and $\varphi \in C^2(\R^d)$ such that $u^* - \varphi$ has a maximum in $x_0$ on $U_\delta(x_0)$ and for $(p,X) \in J^+ u(x_0)$ with $p=D\varphi(x_0)$ and $X \le D^2 \varphi(x_0)$,
\begin{eqnarray*}
\min \left( F(x_0,u^*,p,X, \I_\beta^{1,\delta} [x_0, \varphi(\cdot)] + \I_\beta^{2,\delta} [x_0,p, u^*(\cdot)]  )
, u^* - \mathcal{M}u^* \right) 
&\le & 0 \qquad \mbox{ in } x_0 \in S\\
\min \left( u^* - g, u^* - \mathcal{M}u^* \right) 
&\le & 0 \qquad \mbox{ in } x_0 \in (\R^d \setminus S)
\end{eqnarray*}
A function $u \in \pb$ is a (viscosity) supersolution of (\ref{qvi_elliptic}) if for all $x_0 \in \R^d$ and $\varphi \in C^2(\R^d)$ such that $u_* - \varphi$ has a minimum in $x_0$ on $U_\delta(x_0)$ and for $(q,Y) \in J^- u(x_0)$ with $q=D\varphi(x_0)$ and $Y \ge D^2 \varphi(x_0)$,
\begin{eqnarray*}
\min \left( F(x_0,u_*,q ,Y, \I_\beta^{1,\delta} [x_0, \varphi(\cdot)] + \I_\beta^{2,\delta} [x_0,q, u_*(\cdot)] )
, u_* - \mathcal{M}u_* \right) 
&\ge & 0 \qquad \mbox{ in } x_0 \in S\\
\min \left( u_* - g, u_* - \mathcal{M}u_* \right) 
&\ge & 0 \qquad \mbox{ in } x_0 \in (\R^d \setminus S)
\end{eqnarray*}
A function $u$ is a viscosity solution if it is sub- and supersolution.
\end{definition}

The conditions $p=D\varphi(x_0)$ and $X \le D^2 \varphi(x_0)$ etc. and the maximum condition seem to be superfluous at first view. However, they are needed to ensure consistency of $\varphi$ with the ``local'' derivatives $(p,X)$. 

\begin{proposition}
Definitions 2 and 3 are equivalent.
\end{proposition}

\textbf{Proof: } One direction is obvious. The other direction (see \cite{BaIm07}) uses as vital ingredient that we are considering the \textit{local} maximum. \hfill$\Box$

\subsection{A maximum principle}

Following \cite{BaIm07} we give here a nonlocal theorem which should replace the ``maximum principle''. Prior to this, we have to collect some properties of the intervention operator $\M$ (compare also Lemma \ref{lemma_Mu}):

\begin{lemma}
\label{Lemma_Mu2_etc} 
\begin{enumerate}[(i)]
\item \label{Mu_convex} $\M$ is convex, i.e. for $\lambda \in [0,1]$, $\M ( \lambda a + (1-\lambda) b ) \le \lambda \M a + (1 - \lambda) \M b$
\item \label{Mu_anticonvex} For $\lambda > 0$, $\M ( -\lambda a + (1+\lambda) b ) \ge -\lambda \M a + (1 + \lambda) \M b$ (assuming the latter is not $\infty - \infty$)
\end{enumerate}
\end{lemma}

\textbf{Proof: } Follows easily from $\sup_x (a(x) + b(x)) \le \sup_x a(x) + \sup_x b(x)$ and $\sup_x (a(x) + b(x)) \ge \sup_x a(x) + \inf_x b(x)$, respectively
\hfill$\Box$\\

We need the following nonlocal Jensen-Ishii lemma that can be applied in the PIDE case (compare discussion below):

\begin{lemma}[Lemma 1 in \cite{BaIm07}]
\label{Jensen_Ishii_Lemma}
Let $u\in USC(\R^d)$ and $v \in LSC(\R^d)$, $\varphi \in C^2(\R^{2d})$. If $(x_0,y_0) \in \R^{2d}$ is a zero global maximum point of $u(x) - v(y) - \varphi(x,y)$ and if $p=D_x\varphi(x_0,y_0)$, $q=D_y\varphi(x_0,y_0)$, then for any $K>0$, there exists $\bar{\alpha}(K)>0$ such that, for any $0< \alpha < \bar{\alpha}(K)$, we have: There exist sequences $x_k \to x_0$, $y_k \to y_0$, $p_k \to p$, $q_k \to q$, matrices $X_k, Y_k$ and a sequence of functions $(\varphi_k)$, converging to the function $\varphi_\alpha(x,y) := R^\alpha[\varphi]((x,y),(p,q))$ uniformly in $\R^{2d}$ and in $C^2(B((x_0,y_0),K))$, such that 
\begin{align}
u(x_k)\to u(x_0), v(y_k)\to v(y_0)\\
\label{jensenlemma_eq1} (x_k,y_k) \mbox{ is a global maximum point of } u - v - \varphi_k\\
(p_k,X_k) \in J^+ u(x_k)\\
(q_k,Y_k) \in J^- v(y_k)\\
\label{jensenlemma_maximum_principle_matrix} -\frac{1}{\alpha} I \le \begin{bmatrix} X_k & 0 \\ 0 & -Y_k \end{bmatrix} \le D^2 \varphi_k(x_k,y_k).
\end{align}
Here $p_k = \grad_x \varphi_k(x_k,y_k)$, $q_k = \grad_y \varphi_k(x_k,y_k)$, and $\varphi_\alpha(x_0,y_0) = \varphi(x_0,y_0)$, $\grad \varphi_\alpha(x_0,y_0) = \grad \varphi(x_0,y_0)$. 
\end{lemma}

\begin{rem}
The expression $\varphi_\alpha(x,y) = R^\alpha[\varphi]((x,y),(p,q))$ is the ``modified sup-convolution'' as used by \citet{BaIm07}. For all compacts $C$, $\varphi_\alpha$ converges uniformly to $\varphi$ in $C^2(C)$ as $\alpha \to 0$. This was already used in \cite{BaIm07}, and can be seen by classical arguments using the implicit function theorem. 
\end{rem}

We would obtain a variant of the local Jensen-Ishii lemma (also called maximum principle), if we weren't interested in the sequence $(\varphi_k)$ converging in $C^2$ -- in this case the statement could be expressed in terms of the \textit{limiting} semijets (or ``closures'') $\overline{J}^+$, $\overline{J}^-$ (e.g., $(p,X) \in \overline{J}^+u(x_0)$). However, the local Jensen-Ishii lemma is only useful (in the PDE case), because it can be directly used to deduce, e.g., $F(x_0,u^*(x_0),p,X) \le 0$ by continuity of $F$. Compare also the more detailed explanation in \citet{JaKa06}.

This immediate consequence in the PDE case is a bit more tedious to show in our PIDE case (because the L\'evy measure $\nu$ is possibly singular at $0$), and needs the $C^2$ convergence of the $(\varphi_k)$. The corollary for our impulse control purposes takes the following form:

\begin{corollary}
\label{maximum_principle_impulse_control}
Assume (V\ref{Gamma_cont}), (V\ref{Transact_set_compact}). Let $u$ be a viscosity subsolution and $v$ a viscosity supersolution of (\ref{qvi_elliptic}), and $\varphi \in C^2(\R^{2d})$. If $(x_0,y_0) \in \R^{2d}$ is a global maximum point of $u^*(x)-v_*(y)-\varphi(x,y)$, then, for any $\delta >0$, there exists $\bar{\alpha}$ such that for $0 < \alpha < \bar{\alpha}$, there are $(p,X) \in \overline{J}^+u(x_0)$ and $(q,Y) \in \overline{J}^-v(y_0)$ with
\begin{eqnarray*}
\min \left( F(x_0,u^*(x_0),p,X, \I_\beta^{1,\delta} [x_0, \varphi_\alpha(\cdot,y_0)] + \I_\beta^{2,\delta} [x_0,p, u^*(\cdot)]  )
, u^* - \mathcal{M}u^* \right) 
&\le & 0 \\
\min \left( F(y_0,v_*(y_0),q,Y, \I_\beta^{1,\delta} [y_0, -\varphi_\alpha(x_0,\cdot)] + \I_\beta^{2,\delta} [y_0,q, v_*(\cdot)]  )
, v_* - \mathcal{M}v_* \right) 
&\ge & 0 
\end{eqnarray*}
if $x_0 \in S$ or $y_0 \in S$, respectively. Here $p= \grad_x \varphi(x_0,y_0) = \grad_x \varphi_\alpha(x_0,y_0)$,  $q= - \grad_y \varphi(x_0,y_0) = - \grad_y \varphi_\alpha(x_0,y_0)$, and furthermore,
\begin{equation}
\label{cor_maximum_principle_matrix}
-\frac{1}{\alpha} I \le \begin{bmatrix} X & 0 \\ 0 & -Y \end{bmatrix} \le D^2 \varphi_\alpha(x_0,y_0)= D^2 \varphi(x_0,y_0) + o_\alpha(1).
\end{equation}
\end{corollary}

If none of $x_0$, $y_0$ is in $S$, then we do not need this corollary, because the boundary conditions are not dependent on derivatives of functions $\varphi$.

\begin{rem}
Note that the fact $(p,X) \in \overline{J}^+u(x_0)$ (and the corresponding for the supersolution) is not needed in the statement of the corollary, because the subsolution (supersolution) inequality directly holds by the approximation procedure in the proof. An abstract way of formulating Lemma \ref{Jensen_Ishii_Lemma} and Corollary \ref{maximum_principle_impulse_control} in the style of the local Jensen-Ishii Lemma (only subsolution without impulses) would be to define a new ``limiting superjet'' containing the $(p,X,\varphi_\alpha)$ obtained as limit of the terms in Lemma \ref{Jensen_Ishii_Lemma}. Then Corollary \ref{maximum_principle_impulse_control} could be stated as ``\emph{For $(p,X,\varphi_\alpha)$ in the limiting superjet, $F(x_0,u^*(x_0),p,X,\varphi_\alpha(\cdot),u^*(\cdot)) \le 0$}'' and would follow directly from Lemma \ref{Jensen_Ishii_Lemma}, provided some ``continuity'' of $F: \R^d \times \R \times \R^d \times \spd^d \times C^2 \times \pb \to \R$ hold.
\end{rem}

\textbf{Proof of Corollary \ref{maximum_principle_impulse_control}: } Because of translation invariance, we can assume wlog that $u^*(x_0) - v_*(y_0) - \varphi(x_0,y_0) = 0$. Choose sequences according to Lemma \ref{Jensen_Ishii_Lemma} (applied for $u^*$ and $v_*$), and $K := \max(dist(x_0,U_\delta(x_0)), dist(y_0,U_\delta(y_0))) + 1$. Fix $\alpha \in (0,\bar{\alpha}(K))$.

\textbf{Subsolution case:} If $x_0 \in S$, then $x_k \in S$ for $k$ large, and by Definition 3 and (\ref{jensenlemma_eq1})-(\ref{jensenlemma_maximum_principle_matrix}),
\begin{align}
\label{eq_maxprinc_k_subso_inequality}
\min \left( F(x_k,u^*(x_k),p_k,X_k, \I_\beta^{1,\delta} [x_k, \varphi_k(\cdot,y_k)] + \I_\beta^{2,\delta} [x_k,p_k, u^*(\cdot)]  )
, u^*(x_k) - \mathcal{M}u^*(x_k) \right)
&\le  0 .
\end{align}
First let us prove convergence of the PIDE part $F$. First, $x_k \to x_0$, $u^*(x_k) \to u^*(x_0)$, $p_k \to p$ by Lemma \ref{Jensen_Ishii_Lemma}. $(X_k)$ is contained in a compact set in $R^{d\times d}$ by (\ref{jensenlemma_maximum_principle_matrix}), so it admits a convergent subsequence to an $X$ satisfying (\ref{cor_maximum_principle_matrix}). 
For each $k$, by Proposition \ref{prop_uniqueness_F_continuity} (\ref{prop_F_usc_in_beta}), the supremum in $F$ is attained in a $\beta_k$. Choose another (sub-)subsequence converging to $\beta\in B$.\\
We now only need a reinforced version of (\ref{eq_prop_F_lsc_in_phi_liminf}) in Proposition \ref{prop_uniqueness_F_continuity} for $\I_\beta^{1,\delta}$. By Lebesgue's Dominated Convergence Theorem, (V\ref{Ass_mu_etc_continuous}) and uniform convergence,
\begin{eqnarray*}
&& \lim_{k\to\infty} \int_{|z| < \delta} \varphi_k(x_k + \ell(x_k,\beta_k,z)) - \varphi_k(x_k) - \langle \grad \varphi_k(x_k), \ell(x_k,\beta_k,z) \rangle \, \nu(dz)\\
& = & \int_{|z| < \delta} \varphi_\alpha(x_0 + \ell(x_0,\beta,z)) - \varphi_\alpha(x_0) - \langle \grad \varphi_\alpha(x_0), \ell(x_k,\beta,z) \rangle \, \nu(dz),
\end{eqnarray*}
where the $\nu$-integrable upper estimate can be derived by Taylor expansion and the estimates for $k$ large
\[
\sup_{|z - x|<\kappa_1} |D^2 \varphi_k(z)| \le \sup_{|z - x|<\kappa_1} |D^2 \varphi_\alpha(z)| + \kappa_2
\]
for some $\kappa_1, \kappa_2 > 0$ (recall that $\int_C |z|^2 \nu(dz) < \infty$ for all compacts $C$, and that $U_\delta(y_0) \downarrow 0$ for singular $\nu$). For $\I_\beta^{2,\delta}$, we use Proposition \ref{prop_uniqueness_F_continuity} (\ref{prop_I2_usc}). For the impulse part, we know by Lemma \ref{lemma_Mu} (\ref{Mu_usc}) that $\M u^*$ is usc, so 
\[
\liminf_{k\to\infty} u^*(x_k) - \M u^*(x_k) = u^*(x_0) - \limsup_{k\to\infty} \M u^*(x_k) \ge u^*(x_0) - \M u^*(x_0)
\]
Now we have to combine the estimates derived so far. By iteratively taking subsequences and using (\ref{eq_prop_F_lsc_in_phi_liminf}), we have the desired result for $k \to \infty$ in (\ref{eq_maxprinc_k_subso_inequality}). 

\textbf{Supersolution case:} If $y_0 \in S$, then $y_k \in S$ for $k$ large and by Definition 3 and (\ref{jensenlemma_eq1})-(\ref{jensenlemma_maximum_principle_matrix}),
\begin{align*}
\min \left( F(y_k,v_*(y_k),q_k,Y_k, \I_\beta^{1,\delta} [y_k, -\varphi_k(x_k,\cdot)] + \I_\beta^{2,\delta} [y_k,p_k, v_*(\cdot)]  )
, v_*(y_k) - \M v_*(y_k) \right) 
&\ge  0,
\end{align*}
which means that two separate inequalities for the PIDE part and for the impulse part hold. 
The convergence of the PIDE part is proved in a completely analogous way, except that now (\ref{eq_prop_F_lsc_in_phi_limsup}) in Proposition \ref{prop_uniqueness_F_continuity} is used in a reinforced version (again only needed for $\beta \in B$ fixed). For the impulse part, we know that $\M v_*$ is lsc by Lemma \ref{lemma_Mu} (\ref{Mu_lsc}), so 
\[
\limsup_{k\to\infty} v_*(y_k) - \M v_*(y_k) = v_*(y_0) - \liminf_{k\to\infty} \M v_*(y_k) \le v_*(x_0) - \M v_*(y_0)
\]
\hfill$\Box$

\begin{rem}
By inspecting the proof of Corollary \ref{maximum_principle_impulse_control}, we see that the statement also holds if $u$ and $v$ are subsolution and supersolution, respectively, of different QVIs, provided of course that the conditions are satisfied. This will be used in the proof of the comparison Theorem \ref{viso_comparison}.
\end{rem}


\subsection{A comparison result}
\label{subsec_uniqueness_comparison}

Now we are prepared to give a comparison result (inspired by \cite{Is93}):

\begin{lemma}
\label{lemma_strict_superso}
Assume (V\ref{Gamma_cont}), (V\ref{Transact_set_compact}). 
Let $u$ be a subsolution and $v$ a supersolution of (\ref{qvi_elliptic}), further assume that there is a $w \in \pb \cap C^2(\R^d)$ and a positive function $\kappa: \R^d \to \R$ such that 
\begin{align*}
\min ( - \sup_{\beta \in B} \{ \L^\beta w + f^\beta \}, w - \mathcal{M}w ) &\ge \kappa  \qquad \mbox{in } S\\
\min( w - g, w - \M w ) &\ge \kappa \qquad \mbox{in } \R^d \setminus S.
\end{align*}
Then $v_m := (1 - \frac{1}{m})v + \frac{1}{m} w$ is a supersolution of 
\begin{equation}
\label{qvi_elliptic_strict_superso}
\begin{split}
\min ( - \sup_{\beta \in B} \{ \L^\beta u + f^\beta \}, u - \mathcal{M}u ) - \kappa/m &= 0 \qquad \mbox{in } S\\
\min( u - g, u - \M u ) - \kappa/m &= 0 \qquad \mbox{in } \R^d \setminus S,
\end{split}
\end{equation}
and $u_m := (1 + \frac{1}{m})u - \frac{1}{m} w$ is a subsolution of (\ref{qvi_elliptic}), and of (\ref{qvi_elliptic_strict_superso}) with $-\kappa$ replaced by $+\kappa$.
\end{lemma}

\textbf{Proof: } We use the first viscosity solution definition.\\
First consider the supersolution case. For ease of notation, we write $v$ instead of $v_*$, and so on. Let $\varphi_m \in \pb \cap C^2(\R^d)$, $x_0 \in S$ such that $\varphi_m(x_0) = v_m(x_0)$, $\varphi_m \le v_m$. Choose $\varphi = (\varphi_m - \frac{1}{m}w)(\frac{m}{m-1})$, then $\varphi(x_0) = v(x_0)$ and $\varphi \le v$. We know that $- \sup_{\beta \in B} \{ \L^\beta \varphi + f^\beta \} \ge 0$, and it is sufficient to show that $\L^\beta \varphi_m + f^\beta  \le -\kappa/m$ for all $\beta \in B$, in a point $x_0 \in S$. Using the linearity of $\L$, we obtain
\begin{align*}
0 &\ge  \L^\beta (\varphi_m - \frac{1}{m} w) + \frac{m-1}{m} f^\beta 
= \L^\beta \varphi_m + f^\beta -  \frac{1}{m} (\L^\beta w + f^\beta)
\ge \L^\beta \varphi_m + f^\beta + \frac{\kappa}{m}.
\end{align*}
Because of the convexity of $\M$ (Lemma \ref{Lemma_Mu2_etc} (\ref{Mu_convex})), in any point $x_0 \in \R^d$, 
\begin{align*}
v_m - \M v_m 
&\ge v_m - (1 - \frac{1}{m})\M v - \frac{1}{m} \M w
\ge v_m - (1 - \frac{1}{m})v - \frac{1}{m} \M w\\
& = \frac{1}{m} ( w-  \M w ) > \frac{\kappa}{m}
\end{align*}
It is easy to check that $v_m - g \ge \frac{\kappa}{m}$. 

For the subsolution $u$, the proof proceeds by a case distinction. The reasoning in the impulse part is the same, except that now the anticonvexity of $\M$ (Lemma \ref{Lemma_Mu2_etc} (\ref{Mu_anticonvex})) is used. The PIDE part can be seen (for $\varphi_m(x_0) = u_m(x_0)$, $\varphi_m \ge v_m$, $\varphi = (\varphi_m + \frac{1}{m}w)(\frac{m}{m+1})$) by
\begin{align*}
0 &\le \frac{m+1}{m} \sup_{\beta} [ \L^\beta \varphi + f^\beta ]
= \sup_{\beta} [ \L^\beta \varphi_m + f^\beta +  \frac{1}{m} (\L^\beta w + f^\beta) ]\\
&\le \sup_{\beta} [ \L^\beta \varphi_m + f^\beta ] + \sup_{\beta} [ \frac{1}{m} (\L^\beta w + f^\beta) ]
\le \sup_{\beta} [ \L^\beta \varphi_m + f^\beta ] - \frac{\kappa}{m},
\end{align*}
so we can conclude $- \sup_{\beta} [ \L^\beta \varphi_m + f^\beta ] \le - \frac{\kappa}{m}$. 
\hfill$\Box$

We are going to use the perturbations of sub- and supersolutions to make sure that the maximum of $u_m - v_m$ is attained. So we want to find a $w \ge 0$ growing faster than $|u|$ and $|v|$ as $|x| \to\infty$ (how to find such a $w$ is discussed in Section \ref{sec_model}; the requirements lead to the function $F$ being proper in the sense of \citet{CIL}).

If $\sigma(\cdot,\beta)$, $\mu(\cdot,\beta)$, $f(\cdot,\beta)$, $c$ are Lipschitz continuous, then by classical results (see, e.g., Lemma V.7.1 in \citet{FlSo06}), our function $F$ has the property: 
\begin{quote}
For any $R>0$, there exists a modulus of continuity $\omega_R$, such that, for any $|x|,|y|,|v|\le R$, $l\in\R$ and for any $X,Y \in S^d$ satisfying 
\[
\begin{bmatrix} X & 0 \\ 0 & -Y \end{bmatrix} \le \frac{1}{\varepsilon} \begin{bmatrix} I & -I\\ -I & I \end{bmatrix} + o_\alpha(1)
\]
for some $\varepsilon>0$ ($o_\alpha(1)$ does not depend on $\varepsilon$), then
\[
F(y,v,\varepsilon^{-1}(x-y),Y,l) - F(x,v,\varepsilon^{-1}(x-y),X,l) \le \omega_R(|x-y| + \varepsilon^{-1}|x-y|^2) + o_\alpha(1),
\]
\end{quote}
where $o_\alpha(1)$ again does not depend on $\varepsilon$, and the first term is independent of $\alpha$. In the proof of of Theorem \ref{viso_comparison}, $x_\varepsilon$, $y_\varepsilon$ converge to the same limit $x_0$, so the requirement can be relaxed to locally Lipschitz (as required in Assumption \ref{ass_comparison}) and the property holds only for $x,y$ in a suitable neighbourhood of $x_0$.


Recall that $\pb_p$ is the space of functions in $\pb$ at most polynomially growing with exponent $p$. Assuming essentially that there is a strict supersolution of (\ref{qvi_elliptic}), the following theorem holds:

\begin{theorem}
\label{viso_comparison}
Let Assumptions \ref{ass_existence_uniqueness}, \ref{ass_uniqueness_qvi} and \ref{ass_comparison} be satisfied and $c$ be locally Lipschitz continuous. Assume further that there is a $w\ge 0$ as in Lemma \ref{lemma_strict_superso} (for a constant $\kappa>0$) with $|w(x)|/|x|^p \to \infty$ for $|x| \to\infty$. If $u \in \pb_p(\R^d)$ is a subsolution and $v \in \pb_p(\R^d)$ a supersolution of (\ref{qvi_elliptic}), then $u^* \le v_*$.
\end{theorem}

\begin{corollary}[Viscosity solution: Uniqueness]
Under the same assumptions, there is at most one viscosity solution of (\ref{qvi_elliptic}), and it is continuous.
\end{corollary}

The now following proof of Theorem \ref{viso_comparison} uses the strict sub-/supersolution technique (adapted from \citet{Is93}). We first prove that a maximum can not be attained outside $S$ (because then this would have been because of an impulse back to $S$). Then we use the classical doubling of variables technique, apply the non-local maximum principle, and by a case distinction reduce the problem to a PIDE without impulse part (then adapting techniques of \citet{BaIm07}).\\

\textbf{Proof of Theorem \ref{viso_comparison}: } Write $u$ instead of $u^*$ and $v$ instead of $v_*$ to make the notation more convenient. It is sufficient to prove that $u_m - v_m \le 0$ for all $m$ large (where $u_m, v_m$ are as defined in Lemma \ref{lemma_strict_superso}). Let $m\in\N$ be fixed for the moment. To prove by contradiction, let us assume that $M := \sup_{x \in \R^d} u_m(x) - v_m(x) > 0$. 

\textbf{Step 1. } We want to show that the supremum can not be approximated from within $\R^d \setminus S$. Assume that for each $\varepsilon_1 > 0$, we can find an $\hat{x}=\hat{x}_{\varepsilon_1} \in \R^d \setminus S$ such that $u_m(\hat{x}) - v_m(\hat{x}) + \varepsilon_1 > M$ (and wlog $u_m(\hat{x}) - v_m(\hat{x})>0$). 
By the sub- and supersolution definition, we have
\begin{align*}
\min( u_m(\hat{x}) - g(\hat{x}), u_m(\hat{x}) - \M u_m(\hat{x})) &\le 0\\
\min( v_m(\hat{x}) - g(\hat{x}), v_m(\hat{x}) - \M v_m(\hat{x})) &\ge \kappa/m
\end{align*}
If $u_m(\hat{x}) - g(\hat{x}) \le 0$, then $\kappa/m + u_m(\hat{x}) - v_m(\hat{x}) \le u_m(\hat{x}) - g(\hat{x}) \le 0$ which is already a contradiction. If $u_m(\hat{x}) \le \M u_m(\hat{x})$, then select for $\varepsilon_2 > 0$ a $\hat{\zeta}= \hat{\zeta}_{\varepsilon_1,\varepsilon_2}$ such that $u_m(\Gamma(\hat{x},\hat{\zeta})) + K(\hat{x},\hat{\zeta}) + \varepsilon_2 > \M u_m(\hat{x})$. Then,
\begin{align*}
M - \varepsilon_1 < u_m(\hat{x}) - v_m(\hat{x}) 
&\le u_m(\Gamma(\hat{x},\hat{\zeta})) + K(\hat{x},\hat{\zeta}) + \varepsilon_2 - \kappa/m - K(\hat{x},\hat{\zeta}) - v_m(\Gamma(\hat{x},\hat{\zeta}))\\
&\le \varepsilon_2 - \kappa/m + M,
\end{align*}
which is a contradiction for $\varepsilon_1, \varepsilon_2$ sufficiently small. This shows that the supremum $M$ can not be attained in $\R^d \setminus S$, neither can it be approached from within $\R^d \setminus S$.

\textbf{Step 2. } Now that we are sure we do not have to take into account the boundary conditions, we employ the doubling of variables device as usual. We define for $\varepsilon>0$ and $u_m,v_m$ chosen as in Lemma \ref{lemma_strict_superso}
\[
M_{\varepsilon} = \sup_{x,y \in \R^d} \left( u_m(x) - v_m(y) - \frac{1}{2\varepsilon}  |x-y|^2 \right)
\]
In view of the definition of $w$ and $u_m,v_m$, the maximum is attained in a compact set $C$ (independent of small $\varepsilon$). Choose a point $(x_\varepsilon,y_\varepsilon) \in C$ where the maximum is attained. By applying Lemma 3.1 in \cite{CIL}, we obtain that $\frac{1}{2\varepsilon}  |x_\varepsilon-y_\varepsilon|^2 \to 0$ as $\varepsilon \to 0$, and that $M_\varepsilon \to M = u_m(x_0) - v_m(x_0)$ for all limit points $x_0$ of $(x_\varepsilon)$. We assume from now on wlog that we have chosen a convergent subsequence of $(x_\varepsilon)$, $(y_\varepsilon)$, converging to the same limit $x_0 \in C$. Let $\varepsilon$ small enough such that $x_\varepsilon, y_\varepsilon \in S$ (by Step 1), and that all local estimates in (B\ref{ass_comparison_ints_bounded}), (B\ref{ass_comparison_local_lipschitz}) hold.

Hence, we can apply Corollary \ref{maximum_principle_impulse_control} in $(x_\varepsilon,y_\varepsilon)$ for $\varphi(x,y) = \frac{1}{2\varepsilon}  |x-y|^2$: For any $\delta>0$, there  is a range of $\alpha>0$, for which there are matrices $X,Y$ satisfying (\ref{cor_maximum_principle_matrix}), and $(p,-q) = \grad \varphi(x_\varepsilon,y_\varepsilon)$ (so $p = q = \frac{1}{\varepsilon} (x_\varepsilon - y_\varepsilon)$) such that 
\begin{eqnarray*}
\min \left( F(x_\varepsilon,u_m(x_\varepsilon),p,X, \I_\beta^{1,\delta} [x_\varepsilon, \varphi_\alpha(\cdot,y_\varepsilon)] + \I_\beta^{2,\delta} [x_\varepsilon,p, u_m(\cdot)]  )
, u_m(x_\varepsilon) - \mathcal{M}u_m(x_\varepsilon) \right) 
&\le & 0 \\
\min \left( F(y_\varepsilon,v_m(y_\varepsilon),q,Y, \I_\beta^{1,\delta} [y_\varepsilon, -\varphi_\alpha(x_\varepsilon,\cdot)] + \I_\beta^{2,\delta} [y_\varepsilon,q, v_m(\cdot)]  )
, v_m(y_\varepsilon) - \mathcal{M}v_m(y_\varepsilon) \right) 
&\ge & \frac{\kappa}{m} 
\end{eqnarray*}

\textbf{Case 2a ($u_m(x_\varepsilon) - \mathcal{M}u_m(x_\varepsilon) \le 0$): } Using $v_m(y_\varepsilon) - \mathcal{M}v_m(y_\varepsilon) \ge \frac{\kappa}{m}$, for $\varepsilon>0$ small enough,
\begin{align*}
M 
& = \limsup_{\varepsilon\to 0} ( u_m(x_{\varepsilon}) - v_m(y_{\varepsilon}) )\\
&\le \limsup_{\varepsilon\to 0} \M u_m(x_{\varepsilon}) - \liminf_{\varepsilon\to 0} \M v_m(y_{\varepsilon}) - \frac{\kappa}{m}
\le \M u_m(x_0) - \M v_m(x_0) - \frac{\kappa}{m},
\end{align*}
where we have used the upper and lower semicontinuity of $\M u_m$ and $\M v_m$, respectively (Lemma \ref{Lemma_Mu2_etc}). The contradiction is obtained as in Step 1.

\textbf{Case 2b ($u_m(x_\varepsilon) - \mathcal{M}u_m(x_\varepsilon) > 0$): } It remains to treat the PIDE part
\begin{eqnarray}
\label{Uniqueness_PIDE_subso}
F(x_\varepsilon,u_m(x_\varepsilon),p,X, \I_\beta^{1,\delta} [x_\varepsilon, \varphi_\alpha(\cdot,y_\varepsilon)] + \I_\beta^{2,\delta} [x_\varepsilon,p, u_m(\cdot)]  )
&\le & 0 \\
\label{Uniqueness_PIDE_superso}
F(y_\varepsilon,v_m(y_\varepsilon),q,Y, \I_\beta^{1,\delta} [y_\varepsilon, -\varphi_\alpha(x_\varepsilon,\cdot)] + \I_\beta^{2,\delta} [y_\varepsilon,q, v_m(\cdot)]  )
&\ge & \frac{\kappa}{m}.
\end{eqnarray}
Before we can proceed, we have to compare the integral terms in both inequalities. First note that because $|x + \ell(x,\beta,z) - y|^2 = |x-y|^2 + 2\langle x-y, \ell(x,\beta,z) \rangle + |\ell(x,\beta,z)|^2$, for all $\beta$
\begin{align*}
\I_\beta^{1,\delta} [x_\varepsilon, \varphi(\cdot,y_\varepsilon)] &= \frac{1}{2\varepsilon} \int_{|z|< \delta} | \ell(x_\varepsilon,\beta,z) |^2 \nu(dz) < \infty\\
\I_\beta^{1,\delta} [y_\varepsilon, -\varphi(x_\varepsilon,\cdot)] &= \frac{1}{2\varepsilon} \int_{|z|< \delta} - | \ell(y_\varepsilon,\beta,z) |^2 \nu(dz) < \infty,
\end{align*}
(finite by (V\ref{Ass_ell_integrable}) and definition of $\pb$) so trivially $ \I_\beta^{1,\delta} [x_\varepsilon, \varphi(\cdot,y_\varepsilon)] \le \I_\beta^{1,\delta} [y_\varepsilon, -\varphi(x_\varepsilon,\cdot)] + \frac{1}{\varepsilon} o_\delta(1)$. Because we know that $\varphi_\alpha$ converges to $\varphi$ uniformly in $C^2(C)$ for any compact $C$, we can see analogously to the proof of Corollary \ref{maximum_principle_impulse_control} that $ \I_\beta^{1,\delta} [x_\varepsilon, \varphi_\alpha(\cdot,y_\varepsilon)] \le \I_\beta^{1,\delta} [y_\varepsilon, -\varphi_\alpha(x_\varepsilon,\cdot)] + \frac{1}{\varepsilon} o_\delta(1) + o_\alpha(1)$, where $o_\alpha(1)$ may depend on $\varepsilon$, but is independent of small $\delta$.

Using that $(x_\varepsilon,y_\varepsilon)$ is a maximum point and again $|x + y|^2 = |x|^2 + 2\langle x , y \rangle + |y|^2$,
\begin{equation}
\label{max_ungleichung}
u_m(x_\varepsilon + d ) - u_m(x_\varepsilon) - \frac{1}{\varepsilon} \langle x_\varepsilon - y_\varepsilon , d \rangle
\le v_m(y_\varepsilon + d' ) - v_m(y_\varepsilon) - \frac{1}{\varepsilon} \langle x_\varepsilon - y_\varepsilon , d' \rangle + \frac{1}{2\varepsilon} |d - d'|^2
\end{equation}
where $d, d'$ are arbitrary vectors. We find by integrating (\ref{max_ungleichung}) for all $\beta$ and  $d=\ell(x_\varepsilon,\beta,z)$, $d'=\ell(y_\varepsilon,\beta,z)$ that 
\begin{align*}
\I_\beta^{2,\delta} [x_\varepsilon,p, u_m(\cdot)] \le \I_\beta^{2,\delta} [y_\varepsilon,q, v_m(\cdot)] 
&+ \frac{1}{2\varepsilon} \int_{|z| \ge \delta} |\ell(x_\varepsilon,\beta,z) - \ell(y_\varepsilon,\beta,z)|^2 \,\nu(dz)\\
&+ \int_{|z| \ge 1}  \langle p , \ell(x_\varepsilon,\beta,z) - \ell(y_\varepsilon,\beta,z) \rangle \,\nu(dz).
\end{align*}

We then have by (B\ref{ass_comparison_ints_bounded}) for $\varepsilon>0$ small enough, (denoting $l^1_\beta = \I_\beta^{1,\delta} [x_\varepsilon, \varphi(\cdot,y_\varepsilon)] + \I_\beta^{2,\delta} [x_\varepsilon,p, u_m(\cdot)]$ and $ l^2_\beta = \I_\beta^{1,\delta} [y_\varepsilon, -\varphi(x_\varepsilon,\cdot)] + \I_\beta^{2,\delta} [y_\varepsilon,q, v_m(\cdot)]$) that
\[
l^1_\beta \le l^2_\beta + O(\frac{1}{\varepsilon} |x_\varepsilon- y_\varepsilon|^2) + \frac{1}{\varepsilon} o_\delta(1) + o_\alpha(1)
\]
where $O(\frac{1}{\varepsilon} |x_\varepsilon- y_\varepsilon|^2)$, $o_\delta(1)$ and $o_\alpha(1)$ are independent of $\beta$ because of (B\ref{ass_comparison_ints_bounded}). Likewise, $O(\frac{1}{\varepsilon} |x_\varepsilon- y_\varepsilon|^2)$ is independent of $\delta$ and $\alpha$. Thus
\begin{eqnarray*}
\frac{\kappa}{m}
&\le& F(y_\varepsilon,v_m(y_\varepsilon),q,Y, l^2_\beta ) - F(x_\varepsilon,u_m(x_\varepsilon),p,X, l^1_\beta ) \qquad \mbox{by (\ref{Uniqueness_PIDE_subso}) and (\ref{Uniqueness_PIDE_superso})} \\
&\le& F(y_\varepsilon,v_m(y_\varepsilon),q,Y, l^2_\beta ) - F(x_\varepsilon,v_m(y_\varepsilon),p,X, l^1_\beta ) \qquad \mbox{for small $\varepsilon$ because $c\ge 0$} \\
&\le& F(y_\varepsilon,v_m(y_\varepsilon),q,Y, l^1_\beta ) - F(x_\varepsilon,v_m(y_\varepsilon),p,X, l^1_\beta ) + O(\frac{1}{\varepsilon} |x_\varepsilon- y_\varepsilon|^2) + \frac{1}{\varepsilon} o_\delta(1) + o_\alpha(1)
\end{eqnarray*}
where we have used ellipticity (P\ref{Prop_ellipticity}) and Lipschitz continuity (P\ref{Prop_continuity_l}) in the last component for the $O$ and $o$ values independent of $\beta$. The matrix inequality (\ref{cor_maximum_principle_matrix}) becomes 
\begin{equation}
\label{theo_comparison_matrix_ineq}
-\frac{1}{\alpha} I \le \begin{bmatrix} X & 0 \\ 0 & -Y \end{bmatrix} \le \frac{1}{\varepsilon} \begin{bmatrix} I & -I \\ -I & I \end{bmatrix} + o_\alpha(1).
\end{equation}
By assumption (B\ref{ass_comparison_local_lipschitz}) for $R>0$ large enough ($v_m$ is locally bounded) and $\varepsilon$ small enough,
\[
\frac{\kappa}{m} 
\le \omega_R(|x_{\varepsilon}-y_{\varepsilon}| + {\varepsilon}^{-1}|x_{\varepsilon}-y_{\varepsilon}|^2) + O(\frac{1}{\varepsilon} |x_\varepsilon- y_\varepsilon|^2) + \frac{1}{\varepsilon} o_\delta(1) + o_\alpha(1).
\]
Now let subsequently converge $\delta \to 0$ (because of the special dependence of $\alpha$ -- the smaller $\delta$, the larger $\alpha$ -- this does not affect $\alpha$), and then $\alpha \to 0$.\footnote{By $\alpha\to 0$, we lose the first part of the inequality (\ref{theo_comparison_matrix_ineq}) (so we can not be sure anymore that $X$, $Y$ are bounded because they are dependent on $\alpha$)}. The contradiction is finally obtained by $\varepsilon\to 0$. 
\hfill$\Box$


\subsection{Parabolic case}
\label{subsec_uniqueness_parabolic}

Now let us deduce the parabolic result from the preceding discussion. We will keep the presentation short, only outlining the differences to the elliptic case. We recall the form of the parabolic QVI (where $\partial^+ S_T$ denotes the parabolic nonlocal boundary):
\begin{equation}
\tag{\ref{qvi_parabolic}}
\begin{split}
\min ( - \sup_{\beta \in B} \{ u_t + \L^\beta u + f^\beta \}, u - \M u ) &= 0 \qquad \mbox{in } S_T\\
\min( u - g, u - \M u ) &= 0 \qquad \mbox{in } \partial^+ S_T,
\end{split}
\end{equation}
where for $y = (t,x)$
\begin{multline}
\tag{\ref{L_beta_parabolic}}
\L^\beta u (y) = \frac{1}{2} tr \left(\sigma(y,\beta)\sigma^T(y,\beta)D^2 u(y)\right) + \langle \mu(y,\beta), \grad u(y) \rangle \\
+ \int u(t,x + \ell(y,\beta,z)) - u(y) - \langle \grad u(y), \ell(y,\beta,z) \rangle 1_{|z|< 1} \, \nu(dz).
\end{multline}
The function $F$ then is
\begin{align*}
F(x,r,p,X,(l_\beta)) &= - \sup_{\beta \in B} \left\{ \frac{1}{2} tr \left(\sigma(x,\beta)\sigma^T(x,\beta) X \right) + \langle \mu(x,\beta), p \rangle + f(x,\beta) + l_\beta \right\},
\end{align*}
and the QVI (with the obvious adjustments, compare also section \ref{subsec_existence_parabolic}) then is 
\[
\min ( - u_t(y) + F(y,u(y),\grad_x u(y),D_x^2 u(y), \I_\beta [y, u(t,\cdot)] ), u(y) - \M u(y) ) = 0.
\]

All assumptions and the definition of the space $\pb = \pb([0,T) \times \R^d)$ are as introduced in section \ref{sec_model}. 
The test functions $\varphi$ are now in $C^{1,2}([0,T) \times \R^d)$ (once continuously differentiable in time). Instead of $U_\delta(t_0,x_0) \subset \R^d$, we consider $[0,T) \times U_\delta(t_0,x_0)$, where the $[0,T)$ could in principle be any time interval open in $[0,T)$ containing $t_0$.\\ 
Recall the definition of a viscosity solution of (\ref{qvi_parabolic}) from section \ref{subsec_existence_parabolic}. The original motivation for introducing the different definitions of viscosity solutions was to cater for the singularity in the integral. Because this integral, started in $(t_0,x_0)$, only takes into account values at the time $t_0$, the different definitions in section \ref{subsec_uniqueness_defs} are equivalent in the parabolic case, too (where the time derivative in $t=0$ is only the one-sided derivative).

For the third definition (Definition \ref{viso_def3}), we need the parabolic semijets $P^+u(t,x)$ and $P^-u(t,x)$ on $[0,T)\times \R^d$, e.g.,
\begin{multline*}
P^+u(t,x) = \{ (a,p,X) \in \R \times \R^d \times \spd^d: u(t+s,x+z) \le u(t,x) + a s + \langle p, z \rangle  \\
+ \frac{1}{2} \langle X z, z \rangle + o(|s| + |z|^2) \mbox{ as } s,z \to 0,\; (t+s,z) \in [0,T) \times \R^d \}.
\end{multline*}
The reformulation of Definition \ref{viso_def3} is then ``\textit{For $(a,p,X) \in P^+u(t_0,x_0)$ with $p = D_x \varphi(t_0,x_0)$ and $X \le D^2_x\varphi(t_0,x_0)$, \ldots}'', and in the same way for the supersolution part (we have requirements only on $p$ and $X$ because only they need to be consistent with $\varphi$ as used in the $I^{1,\delta}$ integral).


Finally, we obtain the parabolic maximum principle for impulse control:

\begin{corollary}
\label{maximum_principle_impulse_control_parabolic}
Assume (V\ref{Gamma_cont}), (V\ref{Transact_set_compact}). Let $u$ be a viscosity subsolution and $v$ a viscosity supersolution of (\ref{qvi_parabolic}), and $\varphi \in C^2([0,T) \times \R^{2d})$. If $(t_0,x_0,y_0) \in \R^{2d+1}$ is a global maximum point of $u^*(t,x)-v_*(t,y)-\varphi(t,x,y)$ on $[0,T) \times \R^d$, then, for any $\delta >0$, there exists $\bar{\alpha}$ such that for $0 < \alpha < \bar{\alpha}$, there are $(a,p,X) \in \overline{P}^+u(t_0,x_0)$ and $(b,q,Y) \in \overline{P}^-v(t_0,y_0)$ with
\begin{eqnarray*}
\min \left( -a + F(x_0,u^*(x_0),p,X, \I_\beta^{1,\delta} [x_0, \varphi_\alpha(\cdot,y_0)] + \I_\beta^{2,\delta} [x_0,p, u^*(\cdot)]  )
, u^* - \mathcal{M}u^* \right) 
&\le & 0 \\
\min \left( -b + F(y_0,v_*(y_0),q,Y, \I_\beta^{1,\delta} [y_0, -\varphi_\alpha(x_0,\cdot)] + \I_\beta^{2,\delta} [y_0,q, v_*(\cdot)]  )
, v_* - \mathcal{M}v_* \right) 
&\ge & 0 
\end{eqnarray*}
if $x_0 \in S$ or $y_0 \in S$, respectively. Here $a + b = \varphi_t(t_0,x_0,y_0)$, $p= \grad_x \varphi(t_0,x_0,y_0) = \grad_x \varphi_\alpha(t_0,x_0,y_0)$,  $q= - \grad_y \varphi(t_0,x_0,y_0) = - \grad_y \varphi_\alpha(t_0,x_0,y_0)$, and furthermore,
\begin{equation}
\label{cor_maximum_principle_matrix_parabolic}
-\frac{1}{\alpha} I \le \begin{bmatrix} X & 0 \\ 0 & -Y \end{bmatrix} \le D_{(x,y)}^2 \varphi_\alpha(t_0,x_0,y_0)= D_{(x,y)}^2 \varphi(t_0,x_0,y_0) + o_\alpha(1).
\end{equation}
\end{corollary}

The proof is in no way different, once the parabolic Jensen-Ishii lemma (proved by the same technique as Lemma \ref{Jensen_Ishii_Lemma}) is available. The requirement $a + b = \varphi_t(t_0,x_0,y_0)$ is immediately plausible from the necessary first order criterion in time. Compare also \citet{BaIm07}, \citet{CIL}.

Finally, a comparison theorem can be formulated. Lemma \ref{lemma_strict_superso} is still true in the parabolic case, due to the linearity of the differential operator.

Define $\pb_p = \pb_p([0,T] \times \R^d)$ in the parabolic case by all functions $u \in \pb$, for which there is a (time-independent!) constant $C$ such that $|u(t,x)| \le C(1+|x|^p)$ for all $(t,x) \in [0,T] \times \R^d$. The upper (lower) semicontinuous envelope $u^*$ ($v_*$) is again taken from within $[0,T] \times \R^d$.

\begin{theorem}
\label{viso_comparison_parabolic}
Let Assumptions \ref{ass_existence_uniqueness}, \ref{ass_uniqueness_qvi} and \ref{ass_comparison} be satisfied. Assume further that there is a $w\ge 0$ as in Lemma \ref{lemma_strict_superso} (for a constant $\kappa>0$) with $|w(t,x)|/|x|^p \to \infty$ for $|x| \to\infty$ (uniformly in $t$). If $u \in \pb_p([0,T] \times \R^d)$ is a subsolution and $v \in \pb_p([0,T] \times \R^d)$ a supersolution of (\ref{qvi_parabolic}), then $u^* \le v_*$ on $[0,T] \times \R^d$.
\end{theorem}

\begin{corollary}
Under the same assumptions, there is at most one viscosity solution of (\ref{qvi_parabolic}), and it is continuous on $[0,T] \times \R^d$.
\end{corollary}

\textbf{Proof of Theorem \ref{viso_comparison_parabolic}: } We only point out the differences to the elliptic Theorem \ref{viso_comparison}. Define $M := \sup_{t \in [0,T], x \in \R^d} u_m(t,x) - v_m(t,x)$ and assume it is  $> 0$. Step 1 is proved as in the elliptic case, but on the parabolic boundary $\partial^+ S_T$. For Step 2, we again know that the supremum 
\[
M_\varepsilon = \sup_{t \in [0,T], x,y \in \R^d} \left( u_m(t,x) - v_m(t,y) - \frac{1}{2\varepsilon}  |x-y|^2 \right) 
\]
is attained in a compact set of $[0,T] \times \R^d \times \R^d$ (independent of small $\varepsilon$), say in $(t_\varepsilon,x_\varepsilon,y_\varepsilon)$. For small enough $\varepsilon$, we know by Step 1 that $t_\varepsilon < T$ and $x_\varepsilon, y_\varepsilon \in S$. We proceed as in the elliptic case, and arrive at the PIDE case 2b. All integral estimates hold because $t_\varepsilon$ is fixed at the moment, and the conclusion is exactly the same (the time derivatives $a$ and $b$ cancel out when subtracting the PIDE sub-/supersolution inequalities). The modulus of continuity needs to exist locally uniformly in $t_\varepsilon$ before letting $\varepsilon$ converge to $0$. \hfill$\Box$


\section{Conclusion}
\label{sec_conclusion}

We have shown in the present paper existence and uniqueness of viscosity solutions for impulse control QVI. The results we have obtained are quite general, and the (minimal) assumptions required (basically (local Lipschitz) continuity, continuity of the value function at the boundary, and compactness and ``continuity'' of the transaction set) are sufficient to guarantee a continuous solution on $\R^d$. We note that the Lipschitz continuity assumptions are already needed to ensure existence and uniqueness of the underlying SDE.

The complications to be overcome were mainly 
\begin{itemize}
\item The discontinuous stochastic process and definition of the value function on $\R^d$
\item The possibly singular integral term in the PIDE (arisen from the L\'evy jumps)
\item The additional stochastic control
\end{itemize}

It is our hope that the parabolic and elliptic results presented can be used to great benefit in applications of impulse control without the need to go into details of viscosity solutions (as typically the value function in -- at least financial -- applications will be continuous). The comparison result can also be used to carry out a basic stability analysis for numerical calculations.

Admittedly, our results do not cover all special cases -- but quite frequently, one should be able to extend the results of this paper easily. E.g., state constraints can be handled with a modified framework, where the continuity inside $S$ in general should still hold (see also \citet{LVMP07}). This leaves some room for future research.\\

\textbf{Acknowledgements. } The author wishes to express his thanks to Rüdiger Frey and H. Mete Soner for valuable comments.


\section{Appendix}
\label{sec_appendix}

\begin{lemma}
\label{lemma_taumrho_not0}
Consider a process $X^{\beta^n,t_n,x_n}$ following the SDE (\ref{viso_SDE}), started at $(t_n,x_n)$, and controlled with a stochastic control $\beta^n$. Denote the first exit time for $\rho>0$ by $ \tau_n^\rho := \inf \{ t\ge t_n: | X^{\beta^n,t_n,x_n}_t - x_n | \ge \rho \}$. Suppose that all conditions for existence and uniqueness are satisfied for $X^{\beta^n}$, that (E\ref{ass_ex_higher_variability}) holds and that $(t_n,x_n) \to (t_0,x_0)$. Then for all $\delta>0$ there is a $\hat{\varepsilon}>0$ such that
\[
\limsup_{n\to\infty} \,\P( | \tau_n^\rho - t_0 |< \hat{\varepsilon} ) < \delta.
\]
\end{lemma}

\textbf{Proof: } We want to prove that there is no subsequence in $n$ such that $\P(|\tau_n^\rho - t_n | > \varepsilon) \to 0$ ($n\to \infty$) for all $\varepsilon>0$. Define
\[
K_{n,\varepsilon} := \P(|\tau_n^\rho - t_n| > \varepsilon) = \P ( \sup_{t_n \le s \le t_n + \varepsilon} | X^{\beta^n,t_n,x_n}_s - x_n | < \rho ).
\]
By stochastic continuity, we have for each fixed $n$ that $K_{n,\varepsilon} \to 1$ for $\varepsilon \to 0$. This convergence is uniform in $n$ for the following reasons: By (E\ref{ass_ex_higher_variability}), we can find a constant $\hat{\beta}$ such that $X^{\hat{\beta},t_n,x_n}$ has a higher ``variability'' than $X^{\beta^n,t_n,x_n}$:
\[
\P ( \sup_{t_n \le s \le t_n + \varepsilon} | X^{\hat{\beta},t_n,x_n}_s - x_n | < \rho )
\le \P ( \sup_{t_n \le s \le t_n + \varepsilon} | X^{\beta^n,t_n,x_n}_s - x_n | < \rho )
\]
for all $\varepsilon$ small enough and $n$ large enough. Further, for $n$ large enough such that $|x_n - x_0| < \rho/3$, 
\begin{multline*}
\P ( \sup_{t_n \le s \le t_n + \varepsilon} | X^{\hat{\beta},t_n,x_n}_s - x_n | < \rho ) \\
\ge \P ( \sup_{t_n \le s \le t_n + \varepsilon} | X^{\hat{\beta},t_n,x_n}_s - X^{\hat{\beta},t_0,x_0}_s | < \rho/3 \;\wedge\; \sup_{t_0 \le s \le t_0 + \varepsilon} | X^{\hat{\beta},t_0,x_0}_s - x_0 | < \rho/3 )
\end{multline*}
Finally, by stochastic continuity [in the initial condition] (cf. \citet{GiSk72}, p. 279), the right hand side converges to $1$ for $\varepsilon \to 0$ (uniformly in $n$ for $n$ large enough). This means that for any subsequence $(K_{n_k,\varepsilon})_k$ converging in $k$ for all $\varepsilon>0$, there is an $\hat{\varepsilon}$ such that the limit must lie arbitrarily close to $1$. \hfill$\Box$


\begin{lemma}
\label{lemma_sup_cont}
Let $(X_t)_{t\ge 0}$ be a sequence of random variables, and let $X_t$ converge to $X_0 = x\in \R^d$ in probability for $t\downarrow 0$. Then:
\begin{enumerate}[(i)]
\item \label{lemma_sup_cont_first} For all $\varepsilon>0$, and for all sequences $t_n\to 0$, \[ \P(\sup_{0 \le s \le t_n} |X_s - x| > \varepsilon) \to 0, \qquad n\to\infty\]
\item \label{lemma_sup_cont_second} (\ref{lemma_sup_cont_first}) holds true also for a sequence of positive random variables $(\tau_n)$ converging to $0$ in probability (not necessarily independent of $X$).
\item \label{lemma_sup_cont_uniform} (\ref{lemma_sup_cont_second}) holds true also in the setting and under the assumptions of Lemma \ref{lemma_taumrho_not0}, i.e., 
\[
\P ( \sup_{t_n \le s \le \tau_n} | X^{\beta^n,t_n,x_n}_s - x_n | > \varepsilon ) \to 0
\]
for $(t_n,x_n) \to (t_0,x_0)$, random variables $\tau_n \ge t_n$, $\tau_n \to t_0$ in probability.
\end{enumerate}
\end{lemma}

\textbf{Proof: } (\ref{lemma_sup_cont_first}): By assumption, for all $\varepsilon,\delta>0$ there is a set $U=[0,u]$ for $u>0$ such that for all $s \in U$, $\P(|X_s - x |> \varepsilon)<\delta$. Fix now $\varepsilon, \delta >0$. Dependent on $\omega$, $\gamma>0$ and $t>0$, choose $\hat{t}(\omega,\gamma,t) \in [0,t]$ such that $| X_{\hat{t}(\omega,\gamma,t)} - x | \ge \sup_{0 \le s \le t} |X_t(\omega) - x| - \gamma$.\\
Select now $\gamma:=\varepsilon/2$. Then,
\[
\P( \sup_{0 \le s \le t} | X_s - x| > \varepsilon ) 
\le \P( | X_{\hat{t}(\omega,\gamma,t)} - x | > \varepsilon/2 ) 
\]
Let $U=[0,u]$ be such that $\P( | X_s - x| > \varepsilon/2) < \delta$ for all $s \in U$. Now choose $t:=u$, so $\hat{t}(\omega,\gamma,t) \in [0,u]$, and thus $\P( \sup_{0 \le s \le t} | X_s - x| > \varepsilon ) < \delta$. 

(\ref{lemma_sup_cont_second}): Set for a fixed $\varepsilon>0$ $A_t := \{ \omega: \sup_{0 \le s \le t} |X_s - x| > \varepsilon \}$. We know that $\P(A_t) \to 0$ for $t \to 0$, or, equivalently, $1_{A_t} \to 0$ a.s. Now let $t>0$ be fixed, and $(\tau_n)$ be a sequence of random variables converging to $0$ in probability. Then, a.s., 
\[
1_{A_t} = 1_{A_t} 1_{\{t<\tau_n\}} + 1_{A_t} 1_{\{t \ge \tau_n\}}
\ge 1_{A_t} 1_{\{t<\tau_n\}} + 1_{A_{\tau_n}} 1_{\{t \ge \tau_n\}}
\]
because $1_{A_s} \le 1_{A_t}$ for $s\le t$. The convergence $n\to\infty$ shows that $1_{A_t} \ge \limsup_{n\to\infty} 1_{A_{\tau_n}}$. For $t \to 0$ we get:
\[ 0 \ge \limsup_{n\to\infty} 1_{A_{\tau_n}} \ge 0, \]
thus $1_{A_{\tau_n}} \to 0$ a.s. for $n\to\infty$. 

(\ref{lemma_sup_cont_uniform}): Consider $A_t^n := \{ \omega: \sup_{t_n \le s \le t} |X_s^{\beta_n,t_n,x_n} - x_n| > \varepsilon_1 \}$ in the proof of (\ref{lemma_sup_cont_second}). Then it can be checked in the proof that also $\P(A_{\tau_n}^n) \to 1$ for $n \to \infty$, because by the proof of Lemma \ref{lemma_taumrho_not0}, $\P(A_t^n) \to 1$ for $t \to t_0$, $t \ge t_n$, uniformly in $n$ large.

\hfill$\Box$
\\


\nocite{*}
\bibliography{impulse_viso.bib} 
\bibliographystyle{mynat}

\end{document}